\documentclass[12pt]{article}
\usepackage{amsfonts,amsmath,amsbsy,amssymb,amsthm,graphicx,latexsym,graphics,epsfig,subfigure,color}
\usepackage{mathrsfs}
\usepackage{caption}
\usepackage{subcaption}
\usepackage{array}
\usepackage{booktabs}
\usepackage{tabularx}
\usepackage{longtable,color}
\usepackage{multirow,epsf}
\usepackage{setspace}
\usepackage{epstopdf}
\usepackage{float,wasysym}
\usepackage{color}
\usepackage{enumerate} 
\usepackage{authblk} 

\usepackage{indentfirst}
\usepackage{graphicx}			
\usepackage{mathrsfs}
\usepackage{amsfonts}
\usepackage{amsmath,amssymb,mathrsfs,bm,latexsym,graphicx,color}		
\usepackage{amssymb}
\usepackage[dvipsnames]{xcolor}
\usepackage{natbib} \bibpunct{(}{)}{;}{a}{,}{,} 
\usepackage[normalem]{ulem}
\usepackage[colorlinks=true,urlcolor=blue,citecolor=blue,linkcolor=blue,bookmarks=true]{hyperref} 
\usepackage{titlesec}
\titleformat{\section}[block]{\large\bfseries}{\thesection}{1em}{}
\newcommand{\appendixsection}[1]{
	\refstepcounter{section} 
	\section*{#1} 
	\addcontentsline{toc}{section}{#1} 
	\phantomsection 
}

\newcommand{\appendixref}[1]{Appendix~\ref{#1}}
\usepackage{multicol}
\usepackage{multirow}
\usepackage{tikz,circuitikz} 

\usepackage[T1]{fontenc}
\usepackage[utf8]{inputenc}

\allowdisplaybreaks

\newtheorem{theorem}{Theorem}
\newtheorem{proposition}{Proposition}

\def\Fbar{ {\overline F}}

\def\Hbar{ {\overline H}}

\def\Fbar{ {\overline F}}

\usepackage{tikz,xcolor,hyperref}

\begin{document}
	\title{\bf \Large  Preventive Replacement Policies of Parallel/Series Systems with Dependent Components under Deviation Costs}
\author{	  
	Jiale Niu\thanks{Corresponding author:  E-mail: jiale.niu@outlook.com} \quad   Rongfang Yan\\ 
	{ College of Mathematics and Statistics }\\
	{Northwest Normal University, Lanzhou 730070, China}}


\maketitle

\begin{abstract}
This manuscript studies the preventive replacement policy for a series or parallel system consisting of $n$ independent or dependent heterogeneous components.Firstly, for the age replacement policy, Some sufficient conditions for the existence and uniqueness of the optimal replacement time for both the series and parallel systems are provided. By introducing deviation costs, the expected cost rate of the system is optimized, and the optimal replacement time of the system is extended. Secondly, the periodic replacement policy for series and parallel systems is considered in the dependent case, and a sufficient condition for the existence and uniqueness of the optimal number of periods is provided.  Some numerical examples are given to illustrate and discuss the above preventive replacement policies.
\end{abstract}


\medskip 

\noindent {\bf Mathematicas Subject Classification}: Primary 60E15; Secondary 62G30
\section{Introduction}\label{intro}

 In order to avoid the high costs caused by sudden failures, components or systems must be replaced at planned intervals. The preventive replacement policy involves planning replacement or maintenance before a component or system fails, aiming to minimize sudden failures and reduce maintenance costs. Replacement before and after failure is referred to as preventive replacement (PR) and corrective replacement (CR), respectively. Replacing components too late can cause system failure and high costs, while replacing them too frequently can lead to unnecessary replacement expenses. \cite{barlow1996mathematical} first proposed an age-based replacement policy based on the lifetime of a component or system, which determines the optimal preventive replacement (or repair) time by minimizing the expected cost rate. \cite{berg1976proof} established a preventive maintenance policy for a single-unit system, where replacement occurs either upon failure or a predetermined critical age. It was explained that among all reasonable policies, the age replacement policy itself stands as the optimal criterion for decision-making. \cite{zhao2012optimization} first proposed the preventive replacement policy of "replacement last", developed a new expected cost rate function, and compared it with the classical "replacement first" policy. \cite{zhao2015better}compared the advantages and disadvantages of two preventive replacement models: replacing components at the end of multiple working cycles and replacing components in case of random failures. \cite{levitin2021optimal} studies the modeling and optimization of m-out-of-n standby systems affected by preventive replacement and imperfect component activation. \cite{eryilmaz2023age} studied the age-based preventive replacement policy for any coherent system consisting of mutually independent components with a common discrete life distribution. For comprehensive references on age replacement policy one may refer to \cite{nakagawa2006maintenance,nakagawa2008advanced,zhang2023preventive} and \cite{levitin2024optimizing}.
 
 In most practical situations, it is challenging to meet the assumption of ``replace immediately upon failure'' in the age replacement policy. For example, ocean-going cargo ships must regularly perform navigation tasks and have many critical pieces of equipment and components, such as engines, pumps, navigation systems, and safety equipment. These devices are essential for the ship's normal operation, and their failure could lead to severe consequences, including economic losses and threats to personnel safety. However, due to the limitations of long-distance navigation, if some components suddenly fail, the accompanying personnel usually cannot replace the components. Replacement is typically planned during the ship's port of call to minimize the impact on navigation. In this case, it is more practical to consider a periodic replacement policy. \cite{wireman2004total} has pointed out that the periodic replacement policy is easier to implement in the production system. \cite{nakagawa2008advanced} studied the periodic replacement policy considering the minimum maintenance of components, and optimized the preventive replacement time and the number of working cycles of components by minimizing the expected cost rate. \cite{zhao2014optimal} proposes age and periodic replacement last models with continuous, and discrete policies. \cite{liu2021optimal} studied the optimal periodic preventive maintenance policy for the system affected by external shocks, and explained that the system should be preventively replaced at the optimized periodic time point. In particular, \cite{zhao2024periodic} applied the periodic policy to the backup of database systems and established periodic and random incremental backup policies.
 
 Series and parallel systems, as two very important types of systems in reliability engineering, are widely used in various industrial fields such as power generation
 systems, pump systems, production systems and computing systems. For example, in power systems, transmission lines often operate in series because electricity must pass through each line continuously to be transmitted from the power plant to the user. Therefore, studying the reliability and optimization of series systems is one of the keys to the safe and stable operation of the entire power system(\cite{kundur2007power} and \cite{akhtar2021reliability}). Secondly, in modern computer server systems, parallel system configurations are typically used to ensure high availability and reliability. Multiple components in a parallel system (such as servers, storage devices, and network devices) are redundantly configured. As long as one or more components work properly, the entire system can continue to operate. Parallel design can significantly improve the fault tolerance and reliability of the system(\cite{barroso2022datacenter} and \cite{medara2021energy}). In the past few decades, many scholars have done a lot of work on the reliability optimization (including preventive maintenance) of these two types of systems. \cite{nakagawa2008advanced} aims to minimize the expected cost rate and presents a preventive replacement policy for parallel systems with independent components, focusing on optimizing the number of components in the system and determining the optimal preventive replacement time. Based on the First and Last policy, \cite{nakagawa2012optimization} establishes the optimization of the expected cost rate of the parallel system in terms of the number of components and the replacement time when the number of components is random. \cite{sheu2018generalized} investigated the generalized age maintenance policies for a system with random working times. \cite{wu2017two} discussed the failure process of the system by introducing virtual components. \cite{safaei2020optimal} studied the age replacement policy of repairable series system and parallel system. Based on the minimum expected cost function and the maximum availability function, two optimal age replacement policy are proposed. By considering two types of system failure, \cite{wang2022extended} studied extended preventive replacement models for series and parallel system with n independent non-identical components. \cite{eryilmaz2023optimal} studied the discrete time age replacement policy for a parallel system consisting of components with discrete distributed lifetimes. For further reference on these systems, refer to \cite{xing2021behavior,levitin2023standby} and \cite{xing2024decision}.
 
 The above research on the replacement policy of the system is mostly based on the assumption of component independence. However, in most practical situations, dependencies between components are inevitable. Dependence often occur between components due to shared workloads (such as temperature, humidity, and tasks). For example, the dependent failure was observed in the space development programs in the 1960’s. During the reentry of Gemini spacecraft, one of the two guidance computers failed, and a few minutes later the other one also failed. This is because the temperature inside the two computers was much higher than expected. That is, the two computers failed dependently due to sharing the heat(\cite{ota2017statistical} and \cite{eryilmaz2020optimization}). More unfortunately, on May 31, 2009, shortly after taking off, the flight AF447 crashed into the Atlantic Ocean, killing all 228 people (including 12 crew members) onboard. It is well known that modern aircraft like AF447 are designed with numerous redundant safety systems to ensure absolute safety. In theory, such an accident could only occur if all these safety systems failed simultaneously, which, according to classical reliability theory, is almost impossible. Through the investigation, it was found that the cause of the AF447 accident was erroneous airspeed measurements. For such a failure to occur, all three pitot tubes in the airspeed measurement system would have to fail simultaneously, which is theoretically highly unlikely. However, the reality is that all three pitot tubes froze simultaneously when the aircraft entered a thunderstorm zone. This unfortunate accident also shows that the three pitot tubes failed dependently due to the same load (temperature)(\cite{zeng2023dependent,li2024optimal}). Based on this background, \cite{eryilmaz2020optimization} considered the optimization problem of a system consisting of multiple types of dependent components, and numerically examined how the dependence between the components affects the optimal number of units and replacement time for the system which minimize mean cost rates. Based on two conjecture of \cite{eryilmaz2020optimization}, \cite{torrado2022optimal} presents necessary conditions for the existence of a unique optimal value that minimizes the expected cost rate of two optimization problems for a parallel system consisting of multiple types of components. \cite{xing2020reliability} proposed a combinatorial reliability model for systems undergoing correlated, probabilistic competitions and random failure propagation time for dependent components. \cite{wang2024optimization} used three types of copulas to deal with random maintenance policies for repairable parallel systems with $n$ dependent or independent components. 
 
 To model the dependence between components, the concept of Copula is reviewed below. For a random vector $\bm{X}=(X_1,X_2,\ldots,X_n)$ with the joint distribution function $H$ and respective marginal distribution functions $F_1,F_2, \ldots,F_n$,  the Copula of $X_1,X_2,\ldots,X_n$ is a distribution function $C:[0,1]^n\mapsto [0,1]$, satisfying $$H(\bm{x})=\mathbb{P}(X_1\le x_1,X_2\le x_2,\ldots,X_n\le x_n)=C(F_1(x_1),F_2(x_2),\ldots,F_{n}(x_n)).$$  Similarly, a \emph{survival Copula} of $X_1,X_2,\ldots,X_n$  is a survival function $\hat{C} :[0,1]^n\mapsto [0,1]$, satisfying
 \begin{equation*}
 	\Hbar(\bm{x})=\mathbb{P}(X_1>x_1,X_2>x_2,\ldots,X_n>x_n)=\hat{C}(\Fbar_1(x_1), \Fbar_2(x_2),\ldots,\Fbar_n(x_n)),
 \end{equation*}
 where $ \Hbar(\bm{x})$ is the joint survival function.

 In particular, since Archimedean Copula  has nice mathematical properties,
 which is applied to reliability theory and actuarial science.
 	For a decreasing and continuous function $\phi:[0,1]\mapsto [0,+\infty]$ such that $\phi(0)=+\infty $ and $\phi(1)=0$,  and  denote $\psi=\phi^{-1}$ as the pseudo-inverse of $\phi$. Then
 	$$C_\phi (u_1, u_2,...,u_n)=\psi \Big(\sum_{i=1}^n\phi (u_i)\Big),  \text{ for all } u_i\in [0,1], \quad i=1,2,\dots,n $$
 	is said to be an Archimedean Copula with generator $\phi$ if $(-1)^k\phi^{(k)}(x)\ge 0$ for $k=0,1, \dots,n-2$ and $(-1)^{n-2}\phi^{(n-2)}(x)$ is decreasing and convex.
 	Another family of copulas, widely used in the literature, is the
 	Farlie–Gumbel–Morgenstern (FGM) family which is defined as follows
 	\begin{equation*}
 		C(u_1, u_2,...,u_n)=\prod_{i=1}^{n}u_i\left(1+\theta\prod_{i=1}^{n}(1-u_i) \right), 
 	\end{equation*}
 	for $\theta\in[-1,1]$ with $u_i\in[0,1]$ for $i=1,2,\dots,n$.
    For more discussion on the above Copulas, see \cite{Nelsen2006}.
    
    The distinctiveness and contributions of this research can be highlighted as follows:
    \begin{enumerate}
    	\item [(1)]The age replacement and periodic replacement models of series and parallel systems consisting of dependent heterogeneous components are given(the dependence is modeled by any Copula).
    	\item [(2)]The deviation cost in \cite{zhao2022preventive} is introduced into the replacement policy of series and parallel systems consisting of dependent heterogeneous components, which partially extends the conclusion in \cite{zhao2022preventive}.
    	\item [(3)]Sufficient conditions for the existence and uniqueness of the optimal replacement time is provided for age replacement and periodic replacement models considering deviation costs, with the aim of minimizing the expected cost rate.
    	\item [(4)]Some Copulas satisfying the conditions in this manuscript are given. Numerical examples are given to illustrate and compare the optimal replacement time and the expected cost rate of each policy in this manuscript.
    \end{enumerate}

The remaining part of the paper is organized as follows. Section \ref{s} provides sufficient conditions for the existence and uniqueness of the optimal replacement time to minimize the expected cost rate in the age replacement model of series and parallel systems consisting of dependent heterogeneous components. Section \ref{p} studies the optimal periodic replacement policy for series and parallel systems composed of dependent heterogeneous components. Section \ref{section4} illustrates the conclusions of this work through numerical analysis, and compares the advantages and disadvantages of each preventive replacement policy from the perspectives of dependency, number of components and deviation cost value.  Section \ref{section5} concludes this paper and future directions. And the following Figure \ref{figf1} further illustrates the main work of the manuscript.
\begin{figure}[htbp]
	\centering
	\includegraphics[width=0.8\textwidth]{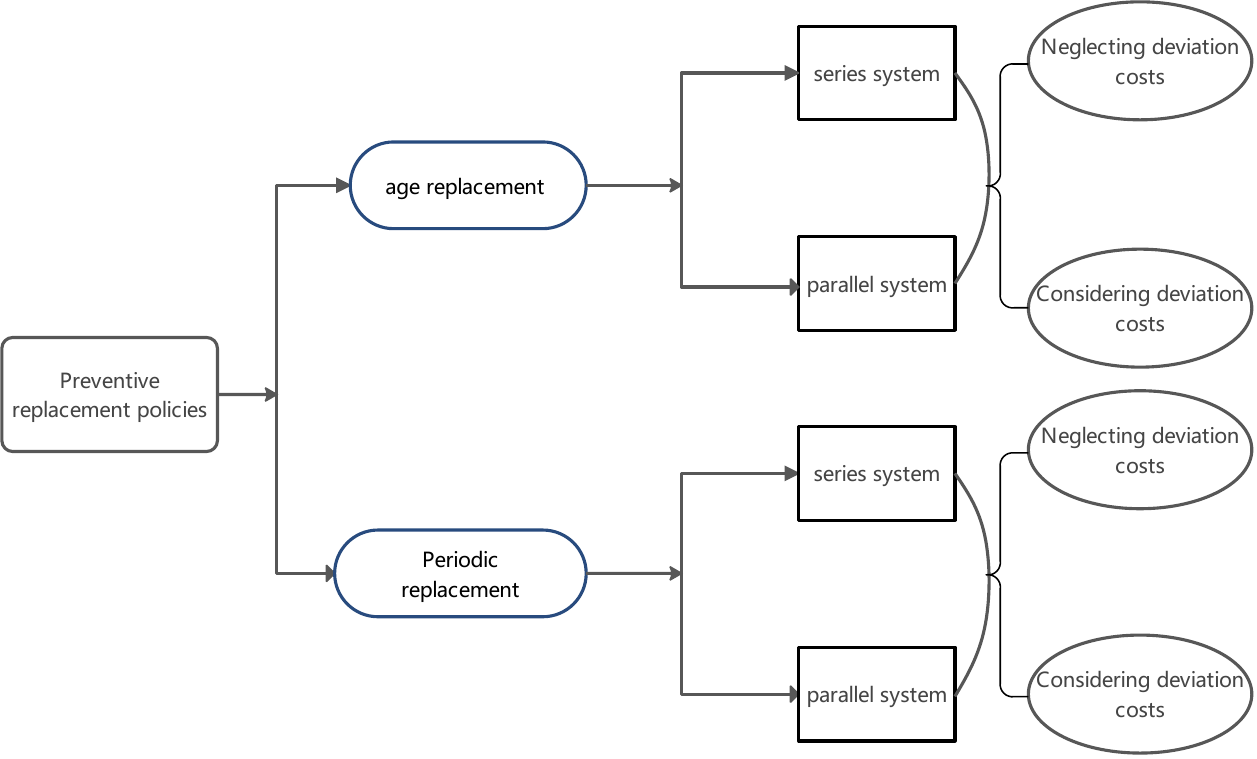}
	\caption{Main research content.}\label{figf1}
\end{figure}
\section{Age-based replacement time T}\label{s}
In this section, we provide the optimal age replacement policy for the series systems and parallel systems considering the deviation cost. By minimizing the expected cost rate, the optimal replacement time is obtained. Assumed that the preventive replacement is carried out at time $T (0< T \leq\infty)$ or correctively at failure time $X$, whichever occurs first. Let $c_{p_i}$ be the preventive replacement cost of the $ith$ component in the system at the planned time $T$, and $c_f$ be a corrective replacement cost for a failed system, where $\sum_{i=1}^{n}c_{p_i}<c_f$. Please see \appendixref{appendix:A} for a more comprehensive proof of this section.

For a series system consisting of $n\geqslant1$ components, the failure of any single component results in the failure of the entire system. Let $X_{1:n}$ denote the lifetime of a series system with n
dependent components assembled by any survival copula
and the heterogeneous lifetimes $X_1,X_2,\dots, X_n$. Then, the reliability function of $X_{1:n}$ can be defined as

\begin{equation*}
	\overline{F}_{1:n}(t)=\hat{C}(\overline{F}_1(t),\overline{F}_2(t),\dots,\overline{F}_n(t)).
\end{equation*}
and the MTTF($\mu_{s}<\infty$) of the system is
\begin{equation*}
	\mu_{s}=\int_{0}^{\infty}\hat{C}(\overline{F}_1(t),\overline{F}_2(t),\dots,\overline{F}_n(t))\mathrm{d}t.
\end{equation*}
Furthermore, the expected cost rate is given by
\begin{eqnarray}\label{(1)}
	C_s(T)&=&\frac{\sum\limits_{i=1}^{n}c_{p_i}\mathbb{P}(X_{1:n}>T)+c_f\mathbb{P}(X_{1:n}\leq T)}{E[\min(X_{1:n},T)]}\nonumber\\
	&=&\frac{c_f-(c_f-\sum\limits_{i=1}^{n}c_{p_i})\hat{C}(\overline{F}_1(T),\overline{F}_2(T),\dots,\overline{F}_n(T))}{\int_{0}^{T}\hat{C}(\overline{F}_1(t),\overline{F}_2(t),\dots,\overline{F}_n(t))\mathrm{d}t}
\end{eqnarray}

The following Theorem \ref{1} shows that conditions on the existence of a unique optimal replacement time $T$ that minimizes the expected cost rate $C_s(T)$ in (\ref{(1)}). To this end, defined the function
\begin{equation*}
	\alpha_i(u_1,u_2,\dots,u_n)=\frac{u_iD_{(i)}H(u_1,u_2,\dots,u_n)}{H(u_1,u_2,\dots,u_n)},
\end{equation*}
where $H(u_1,u_2,\dots,u_n)=\hat{C}(u_1,u_2,\dots,u_n)$ is a generalized domination function such as $\overline{F}_{1:n}(t)=H(\overline{F}_1(t),\overline{F}_2(t),\dots,\overline{F}_n(t))$, and 
\begin{equation*}
	D_{(i)}H(u_1,u_2,\dots,u_n)=\frac{\partial H(u_1,u_2,\dots,u_n)}{\partial u_i},
\end{equation*}
for $i=1,2,\dots,n$.
\begin{theorem}\label{1}
	Let $X_{1:n}$ be a series system with $n$ dependent components assembled by a survival Copula $\hat{C}$ and heterogeneous lifetimes $X_1,X_2,\dots, X_n$. 
	Suppose that the hazard rate $h_{1:n}(\infty)\mu_s>c_f/(c_f-\sum_{i=1}^{n}c_{p_i})$, then there exists a finite and unique optimal replacement time $T^*$ to minimize $C_s(T)$ if $X_i$ is IFR for $i=1,2,\dots,n$ and the function $\alpha_i$ is decreasing in $(0,1)^n$, and the resulting expected cost rate is
	\begin{equation}\label{(3)}
		C_s(T^*)=(c_f-\sum\limits_{i=1}^{n}c_{p_i})h_{1:n}(T^*).
	\end{equation}
\end{theorem}

The definition of deviation time is proposed by \cite{zhao2022preventive}. Assume that the cost of total losses is caused by the deviation time between replacement time $T$ and failure time $X$, as the replacement policy is always planned at earlier or later times than failure times. In the next work, the deviation cost is taken into account in the total cost, and the minimization of the expected cost is studied based on this. In this case, the expected deviation time between $T$ and $X$ is
\begin{eqnarray*}
	&&\int_{0}^{\infty}\mid T-t \mid\mathrm{d}F_{1:n}(t)\\
	&=&\int_{0}^{T}(T-t)\mathrm{d}F_{1:n}(t)+\int_{T}^{\infty}(t-T) \mathrm{d}F_{1:n}(t)\\
	&=&\int_{0}^{T}(1-\hat{C}(\overline{F}_1(t),\!\dots,\!\overline{F}_n(t)))\mathrm{d}t+\int_{T}^{\infty}\hat{C}(\overline{F}_1(t),\!\dots,\!\overline{F}_n(t))\mathrm{d}t
\end{eqnarray*}

Let $c_{d1}$ denotes the deviation downtime cost per unit of time that the system fails before the replacement time $T$, $c_{d2}$ denotes the deviation waste cost per unit of time which caused by the failure of the system later than time $T$, then the total expected cost can be expressed as
\begin{eqnarray*}
	\widetilde{C}_{s1}(T)&=&c_f-(c_f-\sum\limits_{i=1}^{n}c_{p_i})\hat{C}(\overline{F}_1(t),\!\dots,\!\overline{F}_n(t))+	c_{d1}  \int_{0}^{T}( 1-\hat{C}(\overline{F}_1(t),\!\dots,\!\overline{F}_n(t))) \mathrm{d}t\\
	&&+c_{d2}\int_{T}^{\infty}\hat{C}(\overline{F}_1(t),\!\dots,\!\overline{F}_n(t))\mathrm{d}t,
\end{eqnarray*}
and the expected cost rate is
\begin{eqnarray}
	&&\!\!\!{C}_{s1}(T)\nonumber\\
	&=&\!\!\!\frac{c_f\!-\!(c_f\!-\!\sum\limits_{i=1}^{n}c_{p_i})\hat{C}(\overline{F}_1(T),\!\dots,\!\overline{F}_n(T))\!+\!c_{d1}  \int_{0}^{T}( 1-\hat{C}(\overline{F}_1(t),\!\dots,\!\overline{F}_n(t))) \mathrm{d}t\!+\!c_{d2}\mu_s}{\int_{0}^{T}\hat{C}(\overline{F}_1(t),\!\dots,\!\overline{F}_n(t))\mathrm{d}t}-c_{d2}.\nonumber\\
	\label{(5)}
\end{eqnarray}

\begin{theorem}\label{2}
		Let $X_{1:n}$ be a series system with $n$ dependent components assembled by a survival Copula $\hat{C}$ and heterogeneous lifetimes $X_1,X_2,\dots, X_n$. If $X_i$ is IFR for $i=1,2,\dots,n$ and the function $\alpha_i$ is decreasing in $(0,1)^n$, 
 then there exists a finite and unique optimal replacement time $T^*$ to minimize $C_{s1}(T)$, and the resulting expected cost rate is
	\begin{eqnarray}\label{(4)}
		C_{s1}(T^*)=(c_f-\sum\limits_{i=1}^{n}c_{p_i})h_{1:n}(T^*)+c_{d1}\frac{1-\hat{C}(\overline{F}_1(T),\!\dots,\!\overline{F}_n(T))}{\hat{C}(\overline{F}_1(T),\!\dots,\!\overline{F}_n(T))}-c_{d2}.
	\end{eqnarray}
\end{theorem}

The consideration of deviation costs significantly balances the deviation time between replacement and failure in the preventive replacement policy. Additionally, as shown in Theorem \ref{2}, another advantage of incorporating deviation costs is that for exponential distributions with a constant hazard rate, we can also theoretically determine the optimal replacement time to minimize the expected cost rate for the series system. On the other hand, introducing deviation costs eliminates the dependence of the existence and uniqueness condition for the optimal replacement time on $h_{1:n}(\infty)$.

In order to show a practical utility of Theorem \ref{2}, for the ith
component of the series system, suppose the failure time is random variable $X_i$ , and it follows exponential distribution, namely $F_i(t)=P(X_i \leq t)= 1-e^{-\lambda t}~(\lambda >0)$. Moreover, the dependence structure of the system is defined by a Gumbel–Hougaard copula with generator $\phi(t)=(-\ln t)^\theta$ 
for $\theta > 1$. Specially, the Gumbel–Hougaard copula models independence for $\theta=1$. Then, the reliability function and hazard rate function of the corresponding series system, under the above settings, can be expressed as
\begin{equation*}
	\overline{F}_{1:n}(t)=\exp(-(n(-\log \overline{F}(t))^\theta)^{{1}/{\theta}})=(e^{-\lambda t})^{n^{{1}/{\theta}}}
\end{equation*}
and
\begin{equation*}
	h_{1:n}(t)=\lambda n^{{1}/{\theta}}
\end{equation*}
 respectively. The failure rate of the series system is constant under the given conditions. In this scenario, Theorem \ref{1} cannot be used to demonstrate the existence and uniqueness of the system's optimal replacement time. Due to the introduction of deviation costs in Theorem \ref{2}, we can theoretically determine the optimal replacement time for systems with a constant failure rate.

A parallel system consists of $n\geq 1$ components which fails when all units have failed. Let $X_{n:n}$ denote the lifetime of a parallel system with n
dependent components assembled by any copula
and the heterogeneous lifetimes $X_1,X_2,\dots, X_n$. Then, the distribution function of $X_{n:n}$ can be defined as
\begin{equation*}
	{F}_{n:n}(t)={C}({F}_1(t),{F}_2(t),\dots,{F}_n(t)).
\end{equation*}
and the MTTF($\mu_{p}<\infty$) of the system is
\begin{equation*}
	\mu_{p}=\int_{0}^{\infty}(1-{C}({F}_1(t),{F}_2(t),\dots,{F}_n(t)))\mathrm{d}t.
\end{equation*}
Furthermore, the expected cost rate is given by
\begin{eqnarray}\label{(9)}
	C_p(T)&=&\frac{\sum\limits_{i=1}^{n}c_{p_i}\mathbb{P}(X_{n:n}>T)+c_f\mathbb{P}(X_{n:n}\leq T)}{E[\min(X_{n:n},T)]}\nonumber\\
	&=&\frac{\sum\limits_{i=1}^{n}c_{p_i}+(c_f-\sum\limits_{i=1}^{n}c_{p_i}){C}({F}_1(T),\dots,{F}_n(T))}{\int_{0}^{T}(1-{C}({F}_1(t),\dots,{F}_n(t)))\mathrm{d}t}.
\end{eqnarray}
Similarly, to establish the existence and uniqueness condition for the optimal replacement time $T$ that minimizes the expected cost rate $C_p(T)$ in (\ref{(9)}), we define the function
\begin{equation}\label{j1}
	\eta_i(u_1,u_2,\dots,u_n)=\frac{(1-u_i)D_{(i)}H(1-u_1,\dots,1-u_n)}{H(1-u_1,\dots,1-u_n)},
\end{equation}
where $H(1-u_1,\dots,1-u_n)=1-{C}(u_1,u_2,\dots,u_n)$ is a generalized domination function such as $\overline{F}_{n:n}(t)=H(1-{F}_1(t),\dots,1-{F}_n(t))$.

\begin{theorem}\label{3}
Let $X_{n:n}$ be a parallel system with $n$ dependent components assembled by any  Copula ${C}$ and heterogeneous lifetimes $X_1,X_2,\dots, X_n$. 
Suppose that the hazard rate $h_{n:n}(\infty)\mu_{p}>c_f/(c_f-\sum_{i=1}^{n}c_{p_i})$, then there exists a finite and unique optimal replacement time $T^*$ to minimize $C_p(T)$ if $X_i$ is IFR for $i=1,2,\dots,n$ and the function $\eta_i$ is increasing in $(0,1)^n$, and the resulting expected cost rate is
	\begin{equation*}
		C_p(T^*)=(c_f-\sum\limits_{i=1}^{n}c_{p_i})h_{n:n}(T^*).
	\end{equation*}
\end{theorem}

Next, the influence of deviation cost on the optimal replacement time is considered in the preventive age replacement policy of parallel system. The expected cost rate that includes consideration of the deviation cost can be expressed as
\begin{eqnarray}\label{(11)}
	&&\!\!\!{C}_{p1}(T)\nonumber\\
	&=&\!\!\!\frac{\sum\limits_{i=1}^{n}c_{p_i}\!+\!(c_f-\sum\limits_{i=1}^{n}c_{p_i}){C}({F}_1(T),\dots,{F}_n(T))\!+\!	c_{d1}  \int_{0}^{T}{C}({F}_1(t),\dots,{F}_n(t)) \mathrm{d}t+c_{d2}\mu_p}{\int_{0}^{T}(1-{C}({F}_1(t),\dots,{F}_n(t)))\mathrm{d}t}-c_{d2}.\nonumber\\
\end{eqnarray}

\begin{theorem}\label{4}
Let $X_{n:n}$ be a parallel system with $n$ dependent components assembled by any Copula ${C}$ and heterogeneous lifetimes $X_1,X_2,\dots, X_n$. If $X_i$ is IFR for $i=1,2,\dots,n$ and the function $\eta_i$ is increasing in $(0,1)^n$, 
then there exists a finite and unique optimal replacement time $T^*$ to minimize $C_{p1}(T)$, and the resulting expected cost rate is
	\begin{eqnarray*}
		C_{p1}(T^*)=(c_f-\sum\limits_{i=1}^{n}c_{p_i})h_p(T^*)+c_{d1}\frac{{C}({F}_1(T),\dots,{F}_n(T))}{1-{C}({F}_1(T),\dots,{F}_n(T))}-c_{d2}.
	\end{eqnarray*}
\end{theorem}

From now on, for a generalized domination function $H$ defined on $[0,1]^n$, define the sets
\begin{equation*}
	\mathcal{A}=\left\lbrace H: \alpha_i(u_1,u_2,\dots,u_n) ~\text{is decreasing in}~ \mathbf{u}\in [0,1]^n, i=1,2,\dots,n\right\rbrace,
\end{equation*}
and 
\begin{equation*}
	\mathcal{N}=\left\lbrace H: \eta_i(u_1,u_2,\dots,u_n) ~\text{is increasing in}~ \mathbf{u}\in [0,1]^n, i=1,2,\dots,n\right\rbrace.
\end{equation*} 
the non empties of $\mathcal{A}$ and $\mathcal{N}$ have been already showed in \cite{belzunce2001partial} and \cite{navarro2014preservation}.

\begin{proposition}\label{p1}
	\begin{enumerate}
	\item For a series system of $n$ dependent heterogeneous components
	with Archimedean Copula, from Theorem \ref{1}, it is enough to study the monotonicity of
	\begin{equation*}
		\alpha_i(u_1,u_2,\dots,u_n)=\frac{u_i\phi^{'}\left( \sum\limits_{k=1}^{n}\phi^{-1}(u_k)\right) }{\phi^{'}(\phi^{-1}(u_i))\phi\left( \sum\limits_{k=1}^{n}\phi^{-1}(u_k)\right) },
	\end{equation*}
	in ${u_i}\in [0,1]$. For this case, we prove that the Gumbel-Barnett family of copulas for $\theta\in [0,1]$ and the Clayton family of copulas
	for  $\theta\in[-1,0)$ verify such condition;
	\item For FGM copula, the function $\alpha_i(u_1,u_2,\dots,u_n)$ is increasing in $\mathbf{u}\in[0,1]^n$ if and only if $\theta=0$, which corresponds to the case of independence.
  \end{enumerate}
\end{proposition}

For parallel systems consisting of $n$ dependent components, note that Theorem \ref{3} and \ref{4} can be applied to any distributional copula verifying that the function $\eta$ defined in (\ref{j1}) is increasing in $\mathbf{u}\in[0, 1]^n$. In this regard, \cite{torrado2022optimal} prove that the Gumbel–Hougaard family of Copulas for $\theta\geqslant1$ and the Clayton family of Copulas
for $0<\theta\leqslant1$ verify such condition. 
Next, based on \cite{torrado2022optimal}, the proposition \ref{p2} is proposed, which shows that for a parallel system of $n$ dependent components with Clayton Copula, then the function $\eta_i(u_1,u_2,\dots,u_n)$ is increasing in $\mathbf{u}\in[0,1]^n$ if $\theta\in[- 1,0)$.
    \begin{proposition}\label{p2}
		For a parallel system of $n$ dependent components
		with a Clayton Copula, the function $H(u_1,u_2,\dots,u_n)\in\mathcal{N}$ if these $n$ components are homogeneous and $\theta\in[-1,0)$.
	\end{proposition}
	
	In particular, to compare the optimal replacement time and cost of the series system and the parallel system under the same settings, we show that for the Gumbel–Hougaard copulas, when $n$ components of the series system are homogeneous, the corresponding function $\alpha(u)$ is decreasing in $u\in[0,1]$. If the system consists of $n$ homogeneous dependent components with the Gumbel–Hougaard Copula, then the function $\alpha(u)$ can be expressed as
	\begin{equation*}
		\alpha(u)=\frac{n^{1/\theta}u^{n^{1/\theta}-1}u}{u^{n^{1/\theta}}}=n^{1/\theta},
	\end{equation*}
	 i.e. $\alpha(u)$ is decreasing in $u\in[0,1]$.
	 
	 Considering the internal dependence of the system is an important issue in reliability engineering. If the components are originally interdependent, but the wrong assumption of independence is made, the analysis results and expected costs obtained in the replacement policy will differ significantly. This will seriously affect the decision-making of maintenance personnel and managers. For example, we suppose that components $X_1,X_2,\dots,X_8$ follows Weibull distribution, and the distribution function of $X_i$ is 
	 \begin{equation*}
	 	F_i(x)=1-e^{-(\lambda t)^\alpha}, ~\text{for}~t>0,\lambda>0,\alpha>0,i=1,2,\dots,8.
	 \end{equation*}
	 When $c_f=100,c_p=2,c_{d1}=2,c_{d2}=1,\lambda=1.2,\alpha=2.5$, Figure \ref{f1} shows the expected cost rate corresponding to the parallel and series systems in Gumbel-Hougaard Copula, respectively. More detailed discussions of the replacement time and the expected cost rate will be provided in the numerical results analysis section.
\begin{figure}[htbp]
	\centering
	\caption{Comparison of the expected cost rate under different settings.}\label{f1}
	
	\begin{minipage}[b]{0.45\textwidth}
		\centering
		\captionsetup{labelformat=empty}
		\includegraphics[width=\textwidth]{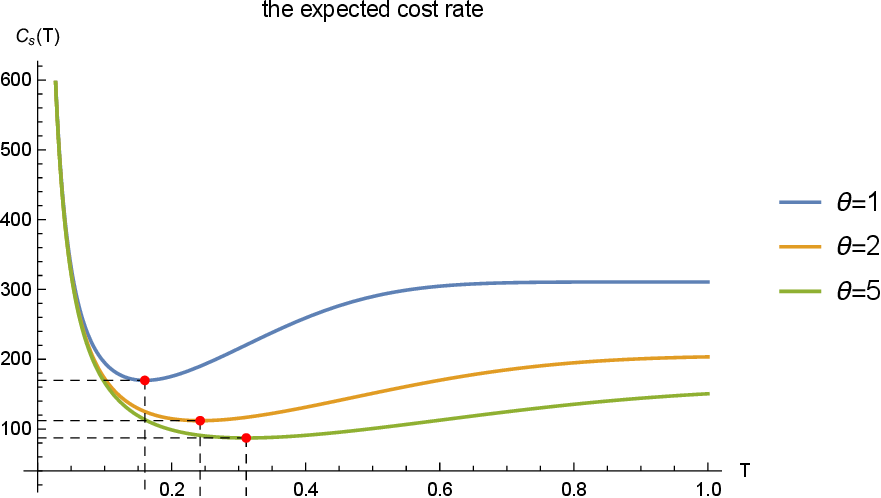}
		\caption{(1):\small{Plot of the expected cost rate $C_s(T)$ in (\ref{(1)}) and dependence parameter $\theta=1,2,5$;}}
	\end{minipage}
	\hfill
	\begin{minipage}[b]{0.45\textwidth}
		\centering
		\captionsetup{labelformat=empty}
		\includegraphics[width=\textwidth]{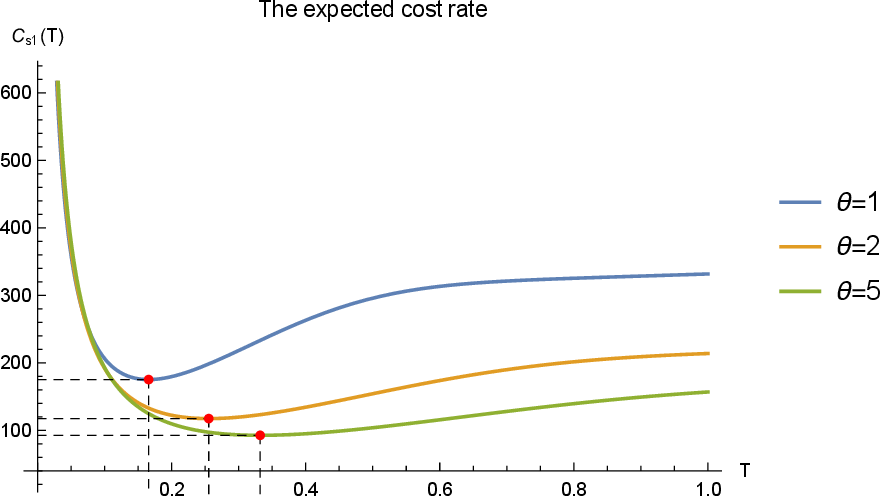}
		\caption{(2):\small{Plot of the expected cost rate $C_{s1}(T)$ in (\ref{(4)}) and dependence parameter $\theta=1,2,5$;}}
	\end{minipage}
	
	\vspace{0.5cm} 
	
	\begin{minipage}[b]{0.45\textwidth}
		\centering
		\captionsetup{labelformat=empty}
		\includegraphics[width=\textwidth]{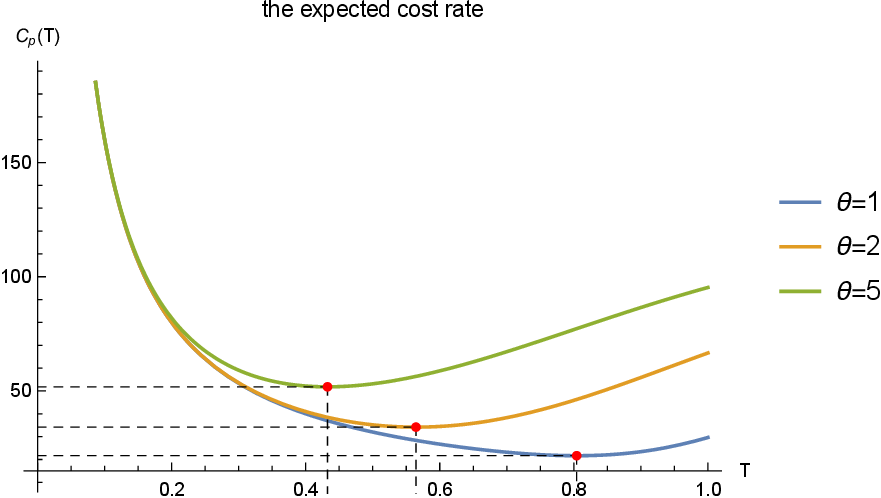}
		\caption{(3):\small{Plot of the expected cost rate $C_p(T)$ in (\ref{(9)}) and dependence parameter $\theta=1,2,5$;}}
	\end{minipage}
	\hfill
	\begin{minipage}[b]{0.45\textwidth}
		\centering
		\captionsetup{labelformat=empty}
		\includegraphics[width=\textwidth]{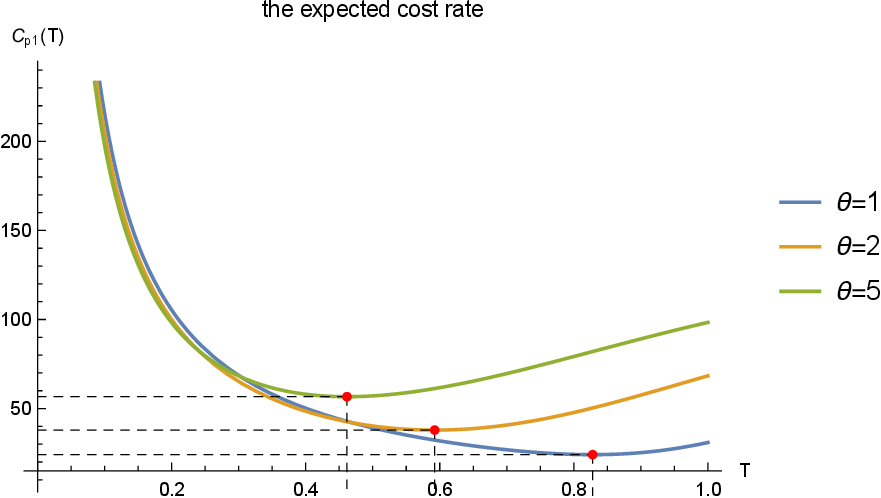}
		\caption{(4):\small{Plot of the expected cost rate $C_{p1}(T)$ in (\ref{(11)}) and dependence parameter $\theta=1,2,5$.}}
	\end{minipage}
	
	\addtocounter{figure}{-4} 
\end{figure}

	It is worth mentioning that the Gumbel-Hougaard copula is positively ordered (\cite{Nelsen2006}), meaning that an increase in the dependence parameter $\theta$ positively reflects the change in dependence among the components.

\section{Periodic replacement time $K\tau$}\label{p}
It has been pointed out that maintenance policies are more easily to be performed at periodic times in production systems (\cite{wireman2004total}). In this section, we plan for a series or parallel system that that preventive replacement is done at periodic times $K\tau~(K=1,2,\dots)$ for a specified
$\tau>0$ or at failure, whichever occurs first. Please see \appendixref{appendix:B} for a more comprehensive proof of this section. Due to the low reliability of the series system, it is necessary to use a periodic replacement policy to check or maintain the system at specified time intervals. Then, for a series system composed of $n$ dependent components, from (\ref{(1)}), the expected cost rate is
\begin{equation}\label{(13)}
	C_s(K)=\frac{c_f-(c_f-\sum\limits_{i=1}^{n}c_{p_i})\hat{C}(\overline{F}_1(K\tau),\dots,\overline{F}_n(K\tau))}{\int_{0}^{K\tau}\hat{C}(\overline{F}_1(t),\dots,\overline{F}_n(t))\mathrm{d}t}.
\end{equation}

%
Next, we find optimum $K^*$ to minimize $C_s(K)$ for given $\tau$.
\begin{theorem}\label{5}
	Let $X_{1:n}$ be a series system with $n$ dependent components assembled by a survival Copula $\hat{C}$ and heterogeneous lifetimes $X_1,X_2,\dots, X_n$. 
	Suppose that the hazard rate $h_{1:n}(\infty)\mu_s>c_f/(c_f-\sum_{i=1}^{n}c_{p_i})$, then there exists a finite and unique $K^*$ to minimize $C_s(K)$ if $X_i$ is IFR for $i=1,2,\dots,n$ and the function $\alpha_i$ is decreasing in $(0,1)^n$.
\end{theorem}

Similarly, due to the periodic replacement policy is implemented at fixed time intervals, it is even more necessary to consider the deviation costs of the system before and after replacement. The following Theorem \ref{6} considers the effect of deviation cost on the expected cost rate of the series system in the periodic replacement policy. From (\ref{(5)}),
\begin{eqnarray}
	&&\!\!\!{C}_{s1}(K)=\nonumber\\
	&&\!\!\!\frac{c_f\!-\!(c_f\!-\!\sum\limits_{i=1}^{n}c_{p_i})\hat{C}(\overline{F}_1(K\tau),\!\dots,\!\overline{F}_n(K\tau))\!+\!c_{d1}  \int_{0}^{K\tau}( 1\!-\!\hat{C}(\overline{F}_1(t),\!\dots,\!\overline{F}_n(t))) \mathrm{d}t\!+\!c_{d2}\mu_s}{\int_{0}^{K\tau}\hat{C}(\overline{F}_1(t),\!\dots,\!\overline{F}_n(t))\mathrm{d}t}\!-\!c_{d2}.\nonumber\\
	\label{(15)}
\end{eqnarray}
\begin{theorem}\label{6}
		Let $X_{1:n}$ be a series system with $n$ dependent components assembled by a survival Copula $\hat{C}$ and heterogeneous lifetimes $X_1,X_2,\dots, X_n$. If $X_i$ is IFR for $i=1,2,\dots,n$ and the function $\alpha_i$ is decreasing in $(0,1)^n$, 
	then there exists a finite and unique minimum $K^*$ to minimize $C_{s1}(K)$ in (\ref{(15)}).
\end{theorem}

For the parallel systems, Theorem \ref{7} determines the optimal periodic replacement time under the periodic replacement policy, with the aim of minimize the expected cost rate.
From (\ref{(9)}), 
\begin{eqnarray}\label{(17)}
	C_p(K)=\frac{\sum\limits_{i=1}^{n}c_{p_i}+(c_f-\sum\limits_{i=1}^{n}c_{p_i}){C}({F}_1(K\tau),\dots,{F}_n(K\tau))}{\int_{0}^{K\tau}(1-{C}({F}_1(t),\dots,{F}_n(t)))\mathrm{d}t}.
\end{eqnarray}
\begin{theorem}\label{7}
	Let $X_{n:n}$ be a parallel system with $n$ dependent components assembled by any Copula ${C}$ and heterogeneous lifetimes $X_1,X_2,\dots, X_n$. Suppose that the hazard rate $h_{n:n}(\infty)\mu_p>c_f/(c_f-\sum_{i=1}^{n}c_{p_i})$, then there exists a finite and unique minimum $K^*$ to minimize $C_s(K)$ if $X_i$ is IFR for $i=1,2,\dots,n$ and the function $\eta_i$ is increasing in $(0,1)^n$.
\end{theorem}

In the following Theorem \ref{8}, the additional deviation cost is considered in the total expected cost of the preventive replacement policy for the parallel system. We will establish a new cost model to find the optimal periodic replacement time that minimizes the expected cost rate. From (\ref{(11)}),
\begin{eqnarray}\label{(19)}
	&&\!\!\!{C}_{p1}(K)\nonumber\\
	&=&\!\!\!\frac{\sum\limits_{i=1}^{n}c_{p_i}\!+\!(c_f\! -\!\sum\limits_{i=1}^{n}c_{p_i}){C}({F}_1(K\tau),\dots,{F}_n(K\tau))\!+\!	c_{d1}  \int_{0}^{K\tau}{C}({F}_1(t),\dots,{F}_n(t)) \mathrm{d}t\!+\!c_{d2}\mu_p}{\int_{0}^{K\tau}(1-{C}({F}_1(t),\dots,{F}_n(t)))\mathrm{d}t}\!-\!c_{d2}.\nonumber\\
\end{eqnarray}
\begin{theorem}\label{8}
	Let $X_{n:n}$ be a parallel system with $n$ dependent components assembled by any Copula ${C}$ and heterogeneous lifetimes $X_1,X_2,\dots, X_n$. If $X_i$ is IFR for $i=1,2,\dots,n$ and the function $\eta_i$ is increasing in $(0,1)^n$, then there exists a finite and unique minimum $K^*$ to minimize $C_{p1}(K)$ in (\ref{(19)}).
\end{theorem}

\begin{figure}[htbp]
	\centering
	\caption{The expected cost rate under different settings.}\label{f2}
	\begin{minipage}[b]{0.45\textwidth}
		\centering
		\captionsetup{labelformat=empty}
		\includegraphics[width=\textwidth]{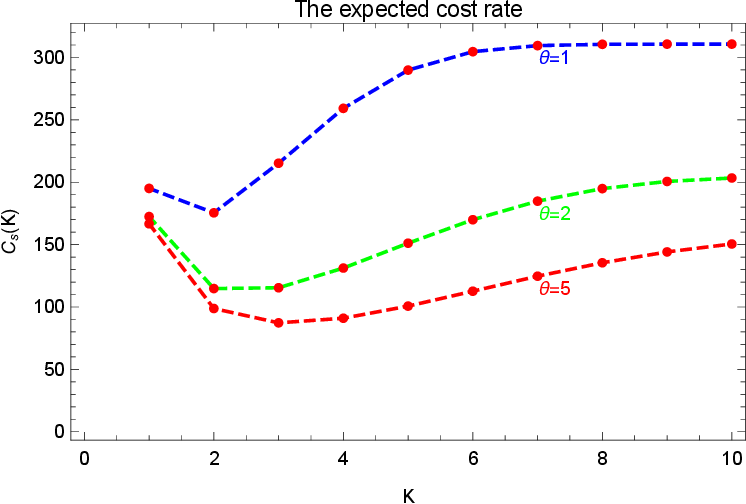}
		\caption{(1):\small{Plot of the expected cost rate $C_s(K)$ in (\ref{(13)}) and dependence parameter $\theta=1,2,5$;}}
		\label{figK:image1}
	\end{minipage}
	\hfill
	\begin{minipage}[b]{0.45\textwidth}
		\centering
		\captionsetup{labelformat=empty}
		\includegraphics[width=\textwidth]{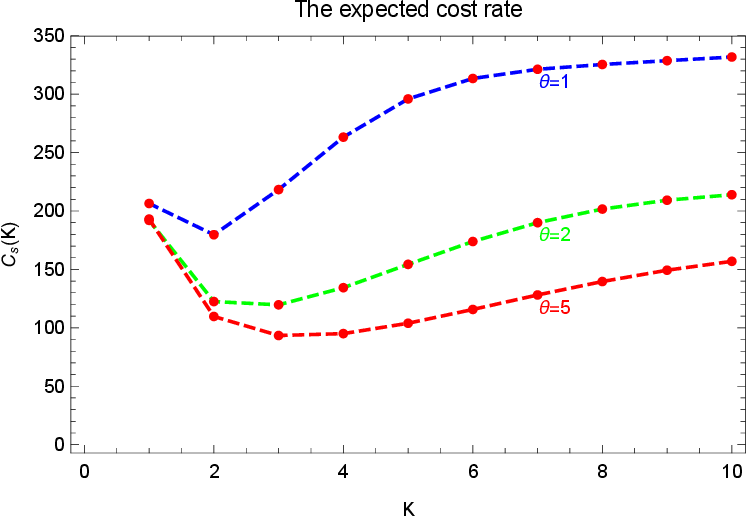}
		\caption{(2):\small{Plot of the expected cost rate $C_{s1}(K)$ in (\ref{(15)}) and dependence parameter $\theta=1,2,5$;}}
		\label{figK:image2}
	\end{minipage}
	
	\vspace{0.5cm} 
	
	\begin{minipage}[b]{0.45\textwidth}
		\centering
		\captionsetup{labelformat=empty}
		\includegraphics[width=\textwidth]{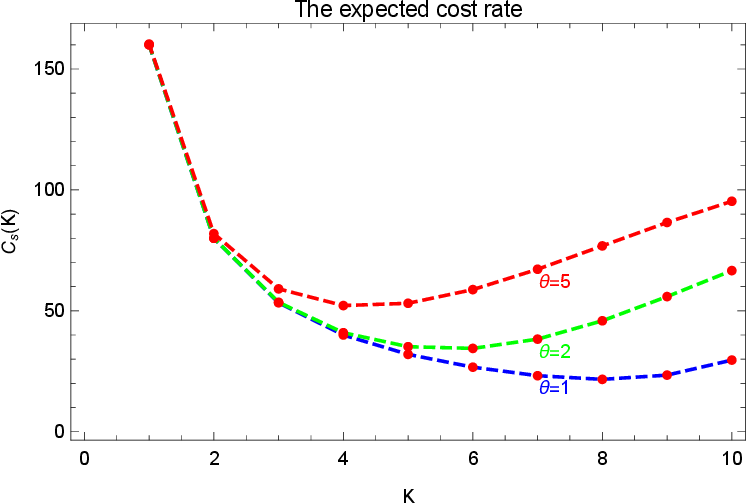}
		\caption{(3):\small{Plot of the expected cost rate $C_p(K)$ in (\ref{(17)}) and dependence parameter $\theta=1,2,5$;}}
		\label{figK:image3}
	\end{minipage}
	\hfill
	\begin{minipage}[b]{0.45\textwidth}
		\centering
		\captionsetup{labelformat=empty}
		\includegraphics[width=\textwidth]{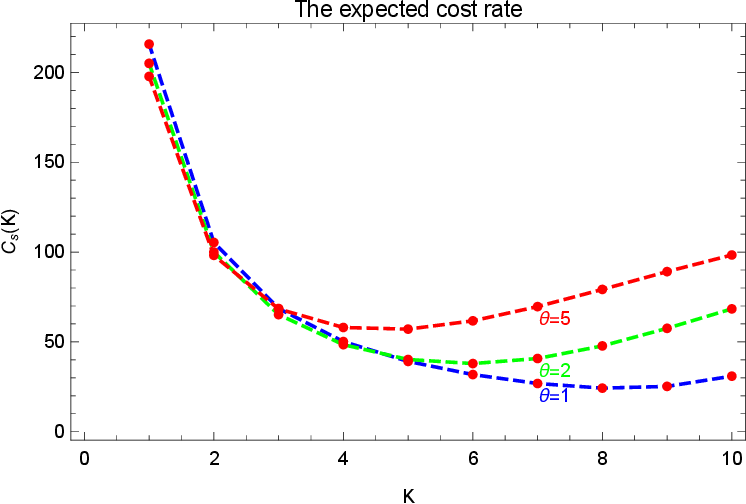}
		\caption{(4):\small{Plot of the expected cost rate $C_{p1}(K)$ in (\ref{(19)}) and dependence parameter $\theta=1,2,5$.}}
		\label{figK:image4}
	\end{minipage}
\end{figure}

To illustrate the shape of the cost function under different dependencies in the periodic replacement policy, we suppose that components $X_1,X_2,\dots,X_8$ follows Weibull distribution. When $c_f=100,c_p=2,c_{d1}=2,c_{d2}=1,\lambda=1.2,\alpha=2.5$ and $\tau=0.1$, Figure \ref{f2} shows the expected cost rate corresponding to the parallel and series systems in Gumbel-Hougaard Copula, respectively.

\section{Numerical results analysis}\label{section4}
In this section, numerical examples are given to illustrate and compare the previous results. Suppose that the failure time $X_i$ of the $i$th component follows the Weibull distribution, that is, $F_i(t)=1- e^{-(\lambda t)^{\alpha}},~i=1,2,\dots,n$. Let $\lambda=0.4,\alpha=2.5$, when $c_f=100,c_{pi}=5$, for the Gumbel-Hougaard copula with $\theta=2$, it is easy to verify that these satisfy the above conditions of Theorems \ref{1}-\ref{4}. Next, the corresponding optimal replacement time and the minimum expected cost rate function in the age replacement policy are presented by Table \ref{b1} and \ref{b2}, respectively, to illustrate the results of Theorem \ref{1},\ref{2},\ref{3},\ref{4}. 
	
\begin{table}[htbp]
	\centering
	\renewcommand{\arraystretch}{1.2} 
	\caption{The optimal $T^*$, $C_s(T^*)$ and $C_{s1}(T^*)$ for $\theta=2$ in the age replacement policy of the series system}
	\label{b1}
	\begin{tabularx}{\textwidth}{c *{6}{>{\centering\arraybackslash}X}}
		\toprule
		\multirow{2}{*}{$n$} & \multicolumn{2}{c}{$c_{d1}=c_{d2}=0$} & \multicolumn{2}{c}{$c_{d1}=2,c_{d2}=1$} & \multicolumn{2}{c}{$c_{d1}=10,c_{d2}=5$}\\
		\cmidrule(lr){2-3} \cmidrule(lr){4-5} \cmidrule(lr){6-7}
		& $\footnotesize T^*$ & $\footnotesize 	C_s(T^*)$ & $\footnotesize T^*$ & $\footnotesize C_{s1}(T^*)$ & $\footnotesize T^*$ & $\footnotesize C_{s1}(T^*)$  \\
		\midrule
		2 & 0.7717 & 21.8277 & 0.8258 & 23.3477 & 0.9927 & 28.3556  \\
		3 & 0.8584 & 29.6197 & 0.8946 & 30.7961 & 1.0127 & 34.9426  \\
		4 & 0.9341 & 36.5439 & 0.9610 & 37.5102 & 1.0491 & 41.0551  \\
		5 & 1.0055 & 42.7784 & 1.0255 & 43.6025 & 1.0926 & 46.7187  \\
		6 & 1.0760 & 48.4142 & 1.0909 & 49.1411 & 1.1406 & 51.9545  \\
		7 & 1.1479 & 53.5065 & 1.1582 & 54.1600 & 1.1921 & 56.7798  \\
		8 & 1.2234 & 58.0923 & 1.2289 & 58.7180 & 1.2469 & 61.2096  \\
		\bottomrule
	\end{tabularx}
\end{table}
	\begin{table}[htbp]
	\centering
	\renewcommand{\arraystretch}{1.2} 
	\caption{The optimal $T^*$, $C_p(T^*)$ and ${C}_{p1}(T^*)$ for $\theta=2$ in the age replacement policy of the parallel system}
	\label{b2}
	\begin{tabularx}{\textwidth}{c *{6}{>{\centering\arraybackslash}X}}
		\toprule
		\multirow{2}{*}{$n$} & \multicolumn{2}{c}{$c_{d1}=c_{d2}=0$} & \multicolumn{2}{c}{$c_{d1}=2,c_{d2}=1$} & \multicolumn{2}{c}{$c_{d1}=10,c_{d2}=5$}\\
		\cmidrule(lr){2-3} \cmidrule(lr){4-5} \cmidrule(lr){6-7}
		& $\footnotesize T^*$ & $\footnotesize 	C_p(T^*)$ & $\footnotesize T^*$ & $\footnotesize C_{p1}(T^*)$ & $\footnotesize T^*$ & $\footnotesize C_{p1}(T^*)$  \\
		\midrule
		2 & 1.0669 & 13.3336 & 1.1391 & 14.6433 & 1.3497 & 18.9166  \\
		3 & 1.3475 & 14.9232 & 1.4038 & 15.8975 & 1.5773 & 19.2541  \\
		4 & 1.5635 & 16.6280 & 1.6099 & 17.4165 & 1.7581 & 20.2180  \\
		5 & 1.7462 & 18.3084 & 1.7823 & 18.9743 & 1.9118 & 21.3905  \\
		6 & 1.8990 & 19.9276 & 1.9338 & 20.5045 & 2.0485 & 22.6307  \\
		7 & 2.0413 & 21.4730 & 2.0721 & 21.9808 & 2.1744 & 23.8775  \\
		8 & 2.1751 & 22.9379 & 2.2025 & 23.3900 & 2.2935 & 25.0990  \\
		\bottomrule
	\end{tabularx}
\end{table}
Based on the numerical results, the following conclusions are drawn.
\begin{enumerate}
	\item [(1)] {By Table \ref{b1}, for a series system, the optimal replacement time $T^*$ decreases with the increase of $n$, which means that the system has shorter lifetime for operations with respect to the increase of $n$. Also, with the increase of $n$, the expected cost rate $C_s(T^*)$ and $C_{s1}(T^*)$ of the system has a significant increase. Therefore, for the series system, it is necessary to reduce the number of components as much as possible without affecting the operation of the system, so as to improve the system reliability and prolong the optimal replacement time. By Table \ref{b2}, for parallel systems, the number of components should be increased as much as possible within the cost limit.}
	\item [(2)] {Incorporating deviation costs into the model will significantly increase $T^*$, $C_s(T^*)(C_p\\(T^*) )$, and $C_{s1}(T^*)(C_{p1}(T^*) )$. 
		This indicates that if the deviation cost $c_{d1}$ and $c_{d2}$ of the system is incorrectly assessed or neglected, it will result in an inaccurate determination of the optimal replacement time $T^*$ and the expected cost rate of the system. 
		In other words, it is necessary to consider the deviation cost in the system replacement policy.}
	\item [(3)]Under the same settings, the parallel system has higher reliability than the series system. Compared with Tables \ref{b1} and \ref{b2}, regardless of the values of $c_{d1}$ and $c_{d2}$, the parallel system always has a longer optimal replacement time and a lower expected cost rate.
\end{enumerate}
In particular, it should be noted that Figure \ref{f1} in Section \ref{s} shows that the dependency between components has a significant impact on the optimal replacement time and the minimum expected cost rate in the replacement policy. In this regard, we assume that the failure time $X_i$ of the $i$th component follows the Weibull distribution, i.e. $F_i(t)=1- e^{-(\lambda t)^{\alpha}},i=1,2,\dots,4$. When $c_f=100,c_{pi}=5$ and $\theta\in\{1,2,4,5,6.5,8.5,15\}$, Tables \ref{b3} and \ref{b4} provide the optimal replacement time and the corresponding expected cost rate for series and parallel systems under different dependent parameters $\theta$.
\begin{table}[htbp]
	\centering
	\renewcommand{\arraystretch}{1.2} 
	\caption{The optimal $T^*$, $C_s(T^*)$ and $C_{s1}(T^*)$ for $\lambda=0.4,\alpha=2.5$ in the age replacement policy of the series system}
	\label{b3}
	\begin{tabularx}{\textwidth}{c *{6}{>{\centering\arraybackslash}X}}
		\toprule
		\multirow{2}{*}{$\theta$} & \multicolumn{2}{c}{$c_{d1}=c_{d2}=0$} & \multicolumn{2}{c}{$c_{d1}=2,c_{d2}=1$} & \multicolumn{2}{c}{$c_{d1}=10,c_{d2}=5$}\\
		\cmidrule(lr){2-3} \cmidrule(lr){4-5} \cmidrule(lr){6-7}
		& $\footnotesize T^*$ & $\footnotesize 	C_s(T^*)$ & $\footnotesize T^*$ & $\footnotesize C_{s1}(T^*)$ & $\footnotesize T^*$ & $\footnotesize C_{s1}(T^*)$  \\
		\midrule
		1 & 0.7080 & 48.2220 & 0.7233 & 49.2067 & 0.7767 & 52.8148  \\
		2 & 0.9341 & 36.5439 & 0.9606 & 37.5302 & 1.0491 & 41.0551  \\
		4 & 1.0730 & 31.8144 & 1.1078 & 32.7993 & 1.2210 & 36.2773  \\
		5 & 1.1032 & 30.9445 & 1.1398 & 31.9294 & 1.2588 & 35.3975  \\
		6.5 & 1.1318 & 30.1627 & 1.1703 & 31.1474 & 1.2948 & 34.6064  \\
		8.5 & 1.1548 & 29.5633 & 1.1948 & 30.5480 & 1.3238 & 33.9996  \\
		15 & 1.1879 & 28.7401 & 1.2302 & 29.7239 & 1.3657 & 33.1652  \\
		\bottomrule
	\end{tabularx}
\end{table}
\begin{table}[htbp]
	\centering
	\renewcommand{\arraystretch}{1.2} 
	\caption{The optimal $T^*$, $C_s(T^*)$ and $C_{s1}(T^*)$ for $\lambda=0.6,\alpha=1.5$ in the age replacement policy of the series system}
	\label{b4}
	\begin{tabularx}{\textwidth}{c *{6}{>{\centering\arraybackslash}X}}
		\toprule
		\multirow{2}{*}{$\theta$} & \multicolumn{2}{c}{$c_{d1}=c_{d2}=0$} & \multicolumn{2}{c}{$c_{d1}=2,c_{d2}=1$} & \multicolumn{2}{c}{$c_{d1}=10,c_{d2}=5$}\\
		\cmidrule(lr){2-3} \cmidrule(lr){4-5} \cmidrule(lr){6-7}
		& $\footnotesize T^*$ & $\footnotesize 	C_s(T^*)$ & $\footnotesize T^*$ & $\footnotesize C_{s1}(T^*)$ & $\footnotesize T^*$ & $\footnotesize C_{s1}(T^*)$  \\
		\midrule
		1 & 0.4468 & 149.1136 & 0.4504 & 150.2213 & 0.4624 & 154.62536 \\
		2 & 0.7092 & 93.9354  & 0.7181 & 95.0426  & 0.7448 & 99.43159  \\
		4 & 0.8935 & 74.5556  & 0.9074 & 75.6631  & 0.9464 & 80.04338  \\
		5 & 0.9358 & 71.1895  & 0.9510 & 72.2961  & 0.9929 & 76.67455  \\
		6.5 & 0.9766 & 68.2170 & 0.9931 & 69.3232 & 1.0380 & 73.69994  \\
		8.5 & 1.0098 & 65.9724 & 1.0274 & 67.0787 & 1.0747 & 71.45403  \\
		15 & 1.0585 & 62.9356 & 1.0777 & 64.0423 & 1.1287 & 68.41568  \\
		\bottomrule
	\end{tabularx}
\end{table}
\begin{table}[htbp]
	\centering
	\renewcommand{\arraystretch}{1.2} 
	\caption{The optimal $T^*$, $C_p(T^*)$ and $C_{p1}(T^*)$ for $\lambda=0.4,\alpha=2.5$ in the age replacement policy of the parallel system}
	\label{b5}
	\begin{tabularx}{\textwidth}{c *{6}{>{\centering\arraybackslash}X}}
		\toprule
		\multirow{2}{*}{$\theta$} & \multicolumn{2}{c}{$c_{d1}=c_{d2}=0$} & \multicolumn{2}{c}{$c_{d1}=2,c_{d2}=1$} & \multicolumn{2}{c}{$c_{d1}=10,c_{d2}=5$}\\
		\cmidrule(lr){2-3} \cmidrule(lr){4-5} \cmidrule(lr){6-7}
		& $\footnotesize T^*$ & $\footnotesize 	C_p(T^*)$ & $\footnotesize T^*$ & $\footnotesize C_{p1}(T^*)$ & $\footnotesize T^*$ & $\footnotesize C_{p1}(T^*)$  \\
		\midrule
		1 & 2.0119 & 11.4888 & 2.0552 & 12.0935 & 2.1908 & 14.2623  \\
		2 & 1.5635 & 16.6280 & 1.6100 & 17.4162 & 1.7581 & 20.2180  \\
		4 & 1.3724 & 21.0679 & 1.4192 & 21.9527 & 1.5691 & 25.0994  \\
		5 & 1.3393 & 22.1888 & 1.3860 & 23.0911 & 1.5356 & 26.3042  \\
		6.5 & 1.3108 & 23.3071 & 1.3573 & 24.2246 & 1.5063 & 27.4971  \\
		8.5 & 1.2899 & 24.2444 & 1.3363 & 25.1730 & 1.4846 & 28.4906  \\
		15 & 1.2629 & 25.6607 & 1.3090 & 26.6036 & 1.4560 & 29.9823  \\
		\bottomrule
	\end{tabularx}
\end{table}
\begin{table}[htbp]
	\centering
	\renewcommand{\arraystretch}{1.2} 
	\caption{The optimal $T^*$, $C_p(T^*)$ and $C_{p1}(T^*)$ for $\lambda=0.6,\alpha=1.5$ in the age replacement policy of the parallel system}
	\label{b6}
	\begin{tabularx}{\textwidth}{c *{6}{>{\centering\arraybackslash}X}}
		\toprule
		\multirow{2}{*}{$\theta$} & \multicolumn{2}{c}{$c_{d1}=c_{d2}=0$} & \multicolumn{2}{c}{$c_{d1}=2,c_{d2}=1$} & \multicolumn{2}{c}{$c_{d1}=10,c_{d2}=5$}\\
		\cmidrule(lr){2-3} \cmidrule(lr){4-5} \cmidrule(lr){6-7}
		& $\footnotesize T^*$ & $\footnotesize 	C_p(T^*)$ & $\footnotesize T^*$ & $\footnotesize C_{p1}(T^*)$ & $\footnotesize T^*$ & $\footnotesize C_{p1}(T^*)$  \\
		\midrule
		1 & 1.3665 & 19.3197 & 1.4112 & 20.2704 & 1.5592 & 23.6470  \\
		2 & 1.0499 & 32.2710 & 1.0918 & 33.3869 & 1.2308 & 37.3880  \\
		4 & 0.9970 & 43.5528 & 1.0360 & 44.6751 & 1.1588 & 48.8264  \\
		5 & 1.0029 & 46.3282 & 1.0406 & 47.4435 & 1.1567 & 51.6144  \\
		6.5 & 1.0154 & 49.0507 & 1.0513 & 50.1587 & 1.1593 & 54.3510  \\
		8.5 & 1.0308 & 51.2922 & 1.0648 & 52.3949 & 1.1646 & 56.6096  \\
		15 & 1.0625 & 54.6011 & 1.0927 & 55.6995 & 1.1773 & 59.9615 \\
		\bottomrule
	\end{tabularx}
\end{table}

Due to the Gumbel-Hougaard family of copulas is positively ordered, i.e. if $C_{\theta_1}\leqslant C_{\theta_2}$ whenever $\theta_1\leqslant\theta_2$. Tables \ref{b3}-\ref{b6} clearly show that changes in the dependency between components significantly impact the optimal replacement time $T^*$ and the corresponding expected average cost rate. This shows that if decision makers and maintenance personnel cannot accurately assess the dependency of the system, it will result in wasted costs due to premature replacement or an increased risk of system failure due to late replacement. But it is noteworthy that the replacement time varies non-monotonically with dependence parameter $\theta$. As shown in Table \ref{b6}, with the increase of $\theta$, the optimal replacement time of the parallel system exhibits a trend of decreasing first and then increasing.

Next, to illustrate the results of Section \ref{p} on periodic replacement policy, assume that the length of a single period is $\tau=0.1$. Suppose that the failure time $X_i$ of the $i$th component follows the Weibull distribution, that is, $F_i(t)=1- e^{-(\lambda t)^{\alpha}},~i=1,2,\dots,n$. When $c_f=100,c_{pi}=5$, assumed that the dependence between components is modeled by the Gumbel-Hougaard copula. it is easy to verify that these satisfy the above conditions of Theorems \ref{5}-\ref{8}. The optimal number of periods $K^*$ and the corresponding expected cost rate $C_s(K^*)(C_p(K^*))$ and ${C}_{s1}(K^*)({C}_{p1}(K^*))$ of the periodic replacement policy of the series and parallel systems are given in Table \ref{b7}-\ref{b10}. Since the periodi replacement policy involves checking and maintaining the system at fixed intervals, it is more convenient for staff to implement compared to the age replacement policy. This means that the periodic replacement policy is more practical and better aligns with real-world scenarios. 
\begin{table}[htbp]
	\centering
	\renewcommand{\arraystretch}{1.2} 
	\caption{The optimal $K^*$, $C_s(K^*)$ and $C_{s1}(K^*)$ for $\lambda=0.4,\alpha=2.5,\theta=2$ in the periodic replacement policy of the series system}
	\label{b7}
	\begin{tabularx}{\textwidth}{c *{6}{>{\centering\arraybackslash}X}}
		\toprule
		\multirow{2}{*}{$n$} & \multicolumn{2}{c}{$c_{d1}=c_{d2}=0$} & \multicolumn{2}{c}{$c_{d1}=2,c_{d2}=1$} & \multicolumn{2}{c}{$c_{d1}=10,c_{d2}=5$}\\
		\cmidrule(lr){2-3} \cmidrule(lr){4-5} \cmidrule(lr){6-7}
		& $\footnotesize K^*$ & $\footnotesize 	C_s(K^*)$ & $\footnotesize K^*$ & $\footnotesize C_{s1}(K^*)$ & $\footnotesize K^*$ & $\footnotesize C_{s1}(K^*)$  \\
		\midrule
		2 & 8 & 21.8480 & 8 & 23.3650& 10 & 28.3569  \\
		3 & 9 & 29.6656 & 9 & 30.7969 & 10 & 34.9470  \\
		4 & 9 & 36.5792 & 10 & 37.5515 & 11 & 41.1240  \\
		5 & 10 & 42.7800 & 10 & 43.6205 & 11 & 46.7202  \\
		6 & 11 & 48.4287 & 11 & 49.1433 & 11 & 52.0014  \\
		7 & 11 & 53.5632 & 12 & 54.2.96 & 12 & 56.7814  \\
		8 & 12 & 58.1040 & 12 & 58.7364 & 12 & 61.2661  \\
		\bottomrule
	\end{tabularx}
\end{table}
\begin{table}[htbp]
	\centering
	\renewcommand{\arraystretch}{1.2} 
	\caption{The optimal $K^*$, $C_p(K^*)$ and ${C}_{p1}(K^*)$ for $\lambda=0.4,\alpha=2.5,\theta=2$ in the periodic replacement policy of the parallel system}
	\label{b8}
	\begin{tabularx}{\textwidth}{c *{6}{>{\centering\arraybackslash}X}}
		\toprule
		\multirow{2}{*}{$n$} & \multicolumn{2}{c}{$c_{d1}=c_{d2}=0$} & \multicolumn{2}{c}{$c_{d1}=2,c_{d2}=1$} & \multicolumn{2}{c}{$c_{d1}=10,c_{d2}=5$}\\
		\cmidrule(lr){2-3} \cmidrule(lr){4-5} \cmidrule(lr){6-7}
		& $\footnotesize K^*$ & $\footnotesize 	C_p(K^*)$ & $\footnotesize K^*$ & $\footnotesize C_{p1}(K^*)$ & $\footnotesize K^*$ & $\footnotesize C_{p1}(K^*)$  \\
		\midrule
		2 & 11 & 13.3477 & 11 & 14.6642 & 14 & 18.9508  \\
		3 & 13 & 14.9489 & 14 & 15.8976 & 16 & 19.2604  \\
		4 & 16 & 16.6414 & 16 & 17.4175 & 18 & 20.2377  \\
		5 & 17 & 18.3259 & 18 & 18.9774 & 19 & 21.3920  \\
		6 & 19 & 19.9276 & 19 & 20.5151 & 20 & 22.6548  \\
		7 & 20 & 21.4876 & 21 & 21.9875 & 22 & 23.8837  \\
		8 & 22 & 22.9428 & 22 & 23.3900 & 23 & 25.0994  \\
		\bottomrule
	\end{tabularx}
\end{table}
\begin{table}[htbp]
	\centering
	\renewcommand{\arraystretch}{1.2} 
	\caption{The optimal $K^*$, $C_s(K^*)$ and ${C}_{s1}(K^*)$ for $\lambda=0.4,\alpha=2.5,n=4$ in the periodic replacement policy of the series system}
	\label{b9}
	\begin{tabularx}{\textwidth}{c *{6}{>{\centering\arraybackslash}X}}
		\toprule
		\multirow{2}{*}{$\theta$} & \multicolumn{2}{c}{$c_{d1}=c_{d2}=0$} & \multicolumn{2}{c}{$c_{d1}=2,c_{d2}=1$} & \multicolumn{2}{c}{$c_{d1}=10,c_{d2}=5$}\\
		\cmidrule(lr){2-3} \cmidrule(lr){4-5} \cmidrule(lr){6-7}
		& $\footnotesize K^*$ & $\footnotesize 	C_s(K^*)$ & $\footnotesize K^*$ & $\footnotesize C_{s1}(K^*)$ & $\footnotesize K^*$ & $\footnotesize C_{s1}(K^*)$  \\
		\midrule
		1 & 7 & 48.2262 & 7  & 49.2271 & 8  & 52.8484  \\
		2 & 9 & 36.5792 & 10 & 37.5515 & 11 & 41.1240  \\
		4 & 11 & 31.8277 & 11 & 32.7777 & 12 & 36.2856  \\
		5 & 11 & 30.9449 & 11 & 31.9337 & 13 & 35.4254  \\
		6.5 & 11 & 30.1791 & 12 & 31.1369 & 13 & 34.6068  \\
		8.5 & 12 & 29.5924 & 12 & 30.5241 & 13 & 34.0082  \\
		15 & 12 & 28.7413 & 12 & 29.7116 & 14 & 33.1809  \\
		\bottomrule
	\end{tabularx}
\end{table}
\begin{table}[htbp]
	\centering
	\renewcommand{\arraystretch}{1.2} 
	\caption{The optimal $K^*$, $C_p(K^*)$ and ${C}_{p1}(K^*)$ for $\lambda=0.4,\alpha=2.5,n=4$ in the periodic replacement policy of the parallel system}
	\label{b10}
	\begin{tabularx}{\textwidth}{c *{6}{>{\centering\arraybackslash}X}}
		\toprule
		\multirow{2}{*}{$\theta$} & \multicolumn{2}{c}{$c_{d1}=c_{d2}=0$} & \multicolumn{2}{c}{$c_{d1}=2,c_{d2}=1$} & \multicolumn{2}{c}{$c_{d1}=10,c_{d2}=5$}\\
		\cmidrule(lr){2-3} \cmidrule(lr){4-5} \cmidrule(lr){6-7}
		& $\footnotesize K^*$ & $\footnotesize 	C_p(K^*)$ & $\footnotesize K^*$ & $\footnotesize C_{p1}(K^*)$ & $\footnotesize K^*$ & $\footnotesize C_{p1}(K^*)$  \\
		\midrule
		1 & 20 & 11.4899 & 21 & 12.1106 & 22 & 14.2632  \\
		2 & 16 & 16.6414 & 16 & 17.4175 & 18 & 20.2377  \\
		4 & 14 & 21.0764 & 14 & 21.9570 & 16 & 25.1107  \\
		5 & 13 & 22.2073 & 14 & 23.0934 & 15 & 26.3199  \\
		6.5 & 13 & 23.3085 & 14 & 24.2455 & 15 & 27.4976  \\
		8.5 & 13 & 24.2456 & 13 & 25.1891 & 15 &  28.4935  \\
		15 & 13 & 25.6768 & 13 & 26.6046 & 15 & 30.0053  \\
		\bottomrule
	\end{tabularx}
\end{table}

\newpage
\section{Concluding remarks}\label{section5}

In this manuscript,  we has systematically examined the age replacement and periodic replacement models for series and parallel systems consisting of dependent heterogeneous components. We consider the age and periodic replacement models for systems with interdependent heterogeneous components, using Copulas to model the dependencies among components. By extending the replacement policy considering deviation costs from \cite{zhao2022preventive} to dependent scenarios, we have expanded the previous results, making the cost assessment for replacement strategies more realistic. Based on the above proposed model, sufficient conditions is provided for the existence and uniqueness of the optimal replacement time to minimize the expected cost rate.

In future research, on the one hand, we plan to consider expanding the models to include multi-objective optimization, balancing cost, reliability, and other performance metrics, will offer a more comprehensive approach to maintenance planning and decision-making. On the other hand, to further study the replacement policy problem, larger and more complex systems, such as coherent systems, series-parallel, and parallel-series systems, are provided to analyze the impact of different maintenance strategies on system performance and lifecycle costs.

\section*{Acknowledgments}
The authors thank two anonymous reviewers and the Editor-in-Chief and  Associate Editor very much for your insightful and constructive comments
and suggestions
, which have greatly improved the presentation of this manuscript.
\section*{Funding}
This research is supported by  the National Natural Science Foundation of China (grant number 11861058). 
\section*{Disclosure statement}
No potential conflict of interest was reported by the author(s).
\appendix
\renewcommand{\thesection}{\Alph{section}} 

\appendixsection{Appendix A}
\label{appendix:A}
Here, we provide the proofs of all the theorems and propositions presented in the section \ref{s}.\newline
\textbf{Proof of Theorem \ref{1}.} 
	Differentiating $C_s(T)$ in (\ref{(1)}) with respect to $T$ and putting it equal to 0, we have
	\begin{eqnarray}
		&&\!\!\!\!\!\!\!\!\!\!\!\!\!\!\!\!\!\!\!\!\!\!\!\!\!\!\!\!\!\!\!\!(c_f-\sum\limits_{i=1}^{n}c_{p_i})h_{1:n}(T)\int_{0}^{T}\hat{C}(\overline{F}_1(t),\overline{F}_2(t),\dots,\overline{F}_n(t))\mathrm{d}t-\Big[ c_f-(c_f-\sum\limits_{i=1}^{n}c_{p_i})\hat{C}(\overline{F}_1(T),\nonumber\\
		&&\!\!\!\!\!\!\!\!\!\!\!\!\!\!\!\!\!\!\!\!\!\!\!\!\!\!\!\!\!\!\!\!\overline{F}_2(T),\dots,\overline{F}_n(T))\Big] =0\nonumber\\
		\text{if and only if}\nonumber\\
		&&\!\!\!\!\!\!\!\!\!\!\!\!\!\!\!\!\!\!\!\!\!\!\!\!\!\!\!\!\!\!\!\!h_{1:n}(T)\int_{0}^{T}\!\!\hat{C}(\overline{F}_1(t),\overline{F}_2(t),\!\dots,\!\overline{F}_n(t))\mathrm{d}t+\hat{C}(\overline{F}_1(T),\overline{F}_2(T),\!\dots,\!\overline{F}_n(T))\!=\!\frac{c_f}{c_f-\sum\limits_{i=1}^{n}c_{p_i}}.\nonumber\\
		\label{(2)}
	\end{eqnarray}
	Then, from Proposition 2.5(i) in \cite{navarro2014preservation}, we can get that $X_{1:n}$ is IFR if $X_i$ is IFR for $i=1,2,\dots,n$ and the function $\alpha_i$ is decreasing in $(0,1)^n$, that is, $h_{1:n}(t)$ is increasing with $t$.
	Next, we show that the left-hand side of (\ref{(2)}) increases with $T$. For convenience, we denote
	\begin{equation*}
		Q_s(T)=h_{1:n}(T)\int_{0}^{T}\!\!\hat{C}(\overline{F}_1(t),\overline{F}_2(t),\!\dots,\!\overline{F}_n(t))\mathrm{d}t+\hat{C}(\overline{F}_1(T),\overline{F}_2(T),\!\dots,\!\overline{F}_n(T)).
	\end{equation*}
	For any $\triangle T>0$,
	\begin{eqnarray*}
		&&\!\!\!Q_s(T+\triangle T)-Q_s(T)\\
		&=&\!\!\!h_{1:n}(T+\triangle T)\int_{0}^{T+\triangle T}\hat{C}(\overline{F}_1(t),\overline{F}_2(t),\!\dots,\!\overline{F}_n(t))\mathrm{d}t+\hat{C}(\overline{F}_1(T+\triangle T),\overline{F}_2(T+\triangle T),\!\dots,\\
		&&\!\!\!\overline{F}_n(T+\triangle T))-h_{1:n}(T)\int_{0}^{T}\hat{C}(\overline{F}_1(t),\overline{F}_2(t),\!\dots,\!\overline{F}_n(t))\mathrm{d}t-\hat{C}(\overline{F}_1(T),\overline{F}_2(T),\!\dots,\!\overline{F}_n(T))\\
		&\geq&\!\!\!h_{1:n}(T\!+\!\triangle T)\int_{0}^{T+\triangle T}\!\!\hat{C}(\overline{F}_1(t),\overline{F}_2(t),\!\dots,\!\overline{F}_n(t))\mathrm{d}t-\!h_{1:n}(T)\!\!\int_{0}^{T}\!\!\hat{C}(\overline{F}_1(t),\overline{F}_2(t),\!\dots,\!\overline{F}_n(t))\mathrm{d}t\\
		&&\!\!\!-h_{1:n}(T\!+\!\triangle T)\int_{T}^{T+\triangle T}\!\!\hat{C}(\overline{F}_1(t),\overline{F}_2(t),\!\dots,\!\overline{F}_n(t))\mathrm{d}t\\
		&=&\!\!\!\left[ h_{1:n}(T+\triangle T)-h_{1:n}(T)\right]\int_{0}^{T}\hat{C}(\overline{F}_1(t),\overline{F}_2(t),\!\dots,\!\overline{F}_n(t))\mathrm{d}t>0, 
	\end{eqnarray*}
	where the first inequality holds due to 
	\begin{eqnarray*}
		\hat{C}(\overline{F}_1(T+\triangle T),\dots,\overline{F}_n(T+\triangle T))-\hat{C}(\overline{F}_1(T),\dots,\overline{F}_n(T))>-h_{1:n}(T+\triangle T),
	\end{eqnarray*}
	caused by the increase of $h_{1:n}(T)$. 
	Additionly, $Q_s(0)=1<c_f/(c_f-\sum_{i=1}^{n}c_{p_i})$, and 
	$Q_s(\infty)=h_s(\infty)\mu_s$.
	Therefore, there exists a finite and unique $T^*~(0<T^*<\infty)$
	that satisfies (\ref{(2)})  and it minimizes $C_s(T)$. Furthermore, combine (\ref{(1)}),  the optimal ACR function is (\ref{(3)}).\newline
\textbf{Proof of Theorem \ref{2}.}
	Differentiating ${C}_{s1}(T)$ in (\ref{(5)}) with respect to $T$ and setting it equal to $0$, 
	\begin{eqnarray}\label{(8)}
		&&\!\!\!\!\!\!\!\!\!\!\!\!\!\!\!\!\!\!\!\!\!\!\!\!\!\!\!\Bigg[ (c_f-\sum\limits_{i=1}^{n}c_{p_i})h_{1:n}(T)+c_{d1}\frac{1-\hat{C}(\overline{F}_1(T),\!\dots,\!\overline{F}_n(T))}{\hat{C}(\overline{F}_1(T),\!\dots,\!\overline{F}_n(T))} \Bigg]\int_{0}^{T}\hat{C}(\overline{F}_1(t),\!\dots,\!\overline{F}_n(t))\mathrm{d}t \nonumber\\
		&&\!\!\!\!\!\!\!\!\!\!\!\!\!\!\!\!\!\!\!\!\!\!\!\!\!\!\!-\Bigg[ c_f-(c_f-\sum\limits_{i=1}^{n}c_{p_i})\hat{C}(\overline{F}_1(T),\!\dots,\!\overline{F}_n(T))+	c_{d1}  \int_{0}^{T}( 1-\hat{C}(\overline{F}_1(t),\!\dots,\!\overline{F}_n(t))) \mathrm{d}t\nonumber\\
		&&\!\!\!\!\!\!\!\!\!\!\!\!\!\!\!\!\!\!\!\!\!\!\!\!\!\!\!+c_{d2}\mu_s\Bigg] =0\nonumber\\
		\text{is equivalent to}\nonumber\\
		&&\!\!\!\!\!\!\!\!\!\!\!\!\!\!\!\!\!\!\!\!\!\!\!\!\!\!\!(c_f-\sum\limits_{i=1}^{n}c_{p_i})Q_s(T)+c_{d1}\varphi_s(T)=c_f+c_{d2}\mu_s, 
	\end{eqnarray}
	where
	\begin{eqnarray*}
		\varphi_s(T)\!\!\!&=&\!\!\!\frac{1\!-\!\hat{C}(\overline{F}_1(T),\!\dots,\!\overline{F}_n(T))}{\hat{C}(\overline{F}_1(T),\!\dots,\!\overline{F}_n(T))}\!\int_{0}^{T}\!\hat{C}(\overline{F}_1(t),\!\dots,\!\overline{F}_n(t))\mathrm{d}t\!-\!\int_{0}^{T}( 1\!-\!\hat{C}(\overline{F}_1(t),\!\dots,\!\overline{F}_n(t))) \mathrm{d}t.
	\end{eqnarray*}
	Hence, 
	\begin{equation*}
		\varphi^\prime_s(T)=\left( \frac{1-\hat{C}(\overline{F}_1(T),\!\dots,\!\overline{F}_n(T))}{\hat{C}(\overline{F}_1(T),\!\dots,\!\overline{F}_n(T))}\right)^ \prime\int_{0}^{T}\hat{C}(\overline{F}_1(t),\!\dots,\!\overline{F}_n(t))\mathrm{d}t>0.
	\end{equation*}
	Then, the left-hand side of (\ref{(8)}) is increases strictly with  T from $0$ to $\infty$, that is, there exists a finite and unique $T^*~(0<T^*<\infty)$ which satisfies (\ref{(8)}), and combine (\ref{(5)}),  and the resulting expected cost rate is
	\begin{equation*}
		C_{s1}(T^*)=(c_f-\sum\limits_{i=1}^{n}c_{p_i})h_{1:n}(T^*)+c_{d1}\frac{1-\hat{C}(\overline{F}_1(T^*),\!\dots,\!\overline{F}_n(T^*))}{\hat{C}(\overline{F}_1(T^*),\!\dots,\!\overline{F}_n(T^*))}-c_{d2}.
	\end{equation*}\newline
\textbf{Proof of Theorem \ref{3}.} 
	Differentiating $C_p(T)$ in (\ref{(9)}) with respect to $T$ and putting it equal to $0$, we have 
	\begin{equation}\label{(10)}
		Q_p(T)=\frac{\sum\limits_{i=1}^{n}c_{p_i}}{c_f-\sum\limits_{i=1}^{n}c_{p_i}}
	\end{equation}
	where 
	\begin{equation*}
		Q_p(T)=h_{n:n}(T)\int_{0}^{T}(1-{C}({F}_1(t),\dots,{F}_n(t)))\mathrm{d}t-{C}({F}_1(t),\dots,{F}_n(t)).
	\end{equation*}
	Then, from Proposition 2.5(i) in \cite{navarro2014preservation}, we can get that $X_{n:n}$ is IFR if $X_i$ is IFR for $i=1,2,\dots,n$ and the function $\eta_i$ is increasing in $(0,1)^n$, that is, $h_{n:n}(t)$ is increasing with $t$.
	For any $\triangle T>0$,
	\begin{eqnarray*}
		&&\!\!\! Q_p(T+\triangle T)-Q_p(T)\\
		&=&\!\!\! h_{n:n}(T\!\!+\!\!\triangle T)\int_{0}^{T+\triangle T}\!(1\!-\!{C}({F}_1(t),\dots,{F}_n(t)))\mathrm{d}t-{C}({F}_1(T+\triangle T),\dots,{F}_n(T+\triangle T))\\
		&&\!\!\!-~h_{n:n}(T)\int_{0}^{T}(1\!-\!{C}({F}_1(t),\dots,{F}_n(t)))\mathrm{d}t+{C}({F}_1(T),\dots,{F}_n(T))\\
		&\geq&\!\!\! h_{n:n}(T\!\!+\!\!\triangle T)\int_{0}^{T+\triangle T}\!(1\!-\!{C}({F}_1(t),\dots,{F}_n(t)))\mathrm{d}t\!-\!h_{n:n}(T)\!\int_{0}^{T}(1\!-\!{C}({F}_1(t),\dots,{F}_n(t)))\mathrm{d}t\\
		&&\!\!\!-~h_{n:n}(T+\triangle T)\int_{T}^{T+\triangle T}(1-{C}({F}_1(t),\dots,{F}_n(t)))\mathrm{d}t\\
		&=&\!\!\!(h_p(T+\triangle T)-h_p(T))\int_{0}^{T}(1-{C}({F}_1(t),\dots,{F}_n(t)))\mathrm{d}t>0.
	\end{eqnarray*}
	Additionly, $Q_p(0)=0$ and $Q_p(\infty)=h_p(\infty)\mu_p-1 $.
	
	Therefore,  if $h_{n:n}(\infty)\mu_{p}>c_f/(c_f-\sum_{i=1}^{n}c_{p_i})$, then there exists a finite and unique $T^*~(0<T^*<\infty)$
	that satisfies (\ref{(10)})  and it minimizes $C_p(T)$. Furthermore, combine (\ref{(9)}),  the optimal ACR function is 
	$C_p(T^*)=(c_f-c_p)h_{n:n}(T^*)$.\newline
\textbf{Proof of Theorem \ref{4}.}
	Differentiating ${C}_{p1}(T)$ in (\ref{(11)}) with respect to $T$ and setting it equal to $0$, 
	\begin{eqnarray}\label{(12)}
		&&\!\!\!\!\!\!\!\!\!\!\!\!\!\!\!\!\!\!\!\!\!\!\!\!\!\!\!\!\!\!\!\!\Bigg[ (c_f-\sum\limits_{i=1}^{n}c_{p_i})h_{n:n}(T)+c_{d1}\frac{{C}({F}_1(T),\dots,{F}_n(T))}{(1-{C}({F}_1(T),\dots,{F}_n(T)))} \Bigg]\int_{0}^{T}(1-{C}({F}_1(t),\dots,\nonumber\\
		&&\!\!\!\!\!\!\!\!\!\!\!\!\!\!\!\!\!\!\!\!\!\!\!\!\!\!\!\!\!\!\!\!{F}_n(t)))\mathrm{d}t -\Bigg[ \sum\limits_{i=1}^{n}c_{p_i}+(c_f-\sum\limits_{i=1}^{n}c_{p_i}){C}({F}_1(T),\dots,{F}_n(T))+	c_{d1}  \int_{0}^{T}{C}({F}_1(t),\dots,\nonumber\\
		&&\!\!\!\!\!\!\!\!\!\!\!\!\!\!\!\!\!\!\!\!\!\!\!\!\!\!\!\!\!\!\!\!{F}_n(t)) \mathrm{d}t+c_{d2}\mu_p\Bigg] =0\nonumber\\
		\text{is equivalent to}\nonumber\\
		&&\!\!\!\!\!\!\!\!\!\!\!\!\!\!\!\!\!\!\!\!\!\!\!\!\!\!\!\!\!\!\!\!(c_f-\sum\limits_{i=1}^{n}c_{p_i})Q_p(T)+c_{d1}\varphi_p(T)=\sum\limits_{i=1}^{n}c_{p_i}+c_{d2}\mu_p, 
	\end{eqnarray}
	where
	\begin{equation*}
		\varphi_p(T)\!=\!\frac{{C}({F}_1(T),\dots,{F}_n(T))}{1\!-\!{C}({F}_1(T),\dots,{F}_n(T))}\!\int_{0}^{T}(1-{C}({F}_1(t),\dots,{F}_n(t)))\mathrm{d}t\!-\!\int_{0}^{T}\!{C}({F}_1(t),\dots,{F}_n(t)) \mathrm{d}t.
	\end{equation*}
	It is easy to get
	\begin{equation*}
		\varphi^\prime_p(T)=\left( \frac{{C}({F}_1(T),\dots,{F}_n(T))}{1-{C}({F}_1(T),\dots,{F}_n(T))}\right)^ \prime\int_{0}^{T}(1-{C}({F}_1(t),\dots,{F}_n(t)))\mathrm{d}t>0.
	\end{equation*}
	Then, the left-hand side of (\ref{(12)}) is increases strictly with  $T$ from $0$ to $\infty$, that is, there exists a finite and unique $T^*~(0<T^*<\infty)$ which satisfies (\ref{(12)}), and combine (\ref{(11)}),  and the resulting cost rate is
	\begin{equation*}
		C_{p1}(T^*)=(c_f-\sum\limits_{i=1}^{n}c_{p_i})h_{n:n}(T^*)+c_{d1}\frac{{C}({F}_1(T),\dots,{F}_n(T))}{1-{C}({F}_1(T),\dots,{F}_n(T))}-c_{d2}.
	\end{equation*}
\textbf{Proof of Proposition \ref{p1}.}
	\begin{enumerate}
		\item First, we consider the Gumbel-Barnett copula family; for its proof, please refer to \cite{yan2022component}.
		
		Secondly, for the Clayton family of copulas, 
		\begin{equation*}
			H(u_1,u_2,\dots, u_n)=\max \left[ \left(u_1^{-\theta}+u_2^{-\theta}+\dots +u_n^{-\theta}-n+1,0 \right)\right] ^{-1/\theta },~~for~u_i\in[0,1].
		\end{equation*}
		To show that $\alpha_i(u_1,u_2,\dots,u_n)$ is decreasing in $\mathbf{u}\in [0,1]^n, ~i=1,2,\dots,n$, it suffices to consider the case where $u_1^{-\theta}+u_2^{-\theta}+\dots +u_n^{-\theta}-n+1>0$. For any $i= {1,2,\dots,n}$,
		\begin{eqnarray*}
			\alpha_i(u_1,u_2,\dots,u_n)&=&\frac{u_iD_{(i)}H(u_1,u_2,\dots,u_n)}{H(u_1,u_2,\dots,u_n)}\\
			&=&\frac{\left(u_1^{-\theta}+u_2^{-\theta}+\dots +u_n^{-\theta}-n+1 \right)^{-1/\theta-1}u_i^{-\theta}}{\left(u_1^{-\theta}+u_2^{-\theta}+\dots +u_n^{-\theta}-n+1 \right)^{-1/\theta }}\\
			&=&\left(u_1^{-\theta}+u_2^{-\theta}+\dots +u_n^{-\theta}-n+1 \right)^{-1}u_i^{-\theta},~~\theta \in[-1,0).
		\end{eqnarray*}
		Taking the partial derivative of $\alpha_i(u_1,u_2,\dots,u_n)$ with respect to $u_i,$ we have
		\begin{eqnarray*}
			&&\!\!\!\frac{\partial \alpha_i(u_1,u_2,\dots,u_n)}{\partial u_i}\\
			&=&\!\!\!\theta(u_1^{-\theta}+\dots +u_n^{-\theta}-n+1)^{-2} u_i^{-2\theta-1}-\theta (u_1^{-\theta}+\dots +u_n^{-\theta}-n+1)^{-1}u_i^{-\theta-1}\\
			&=&\!\!\!-\theta(u_1^{-\theta}+\dots +u_n^{-\theta}-n+1)^{-1}u_i^{-\theta-1}\left[ 1-(u_1^{-\theta}+\dots +u_n^{-\theta}-n+1)^{-1}u_1^{-\theta}\right] \\
			&\overset{\rm sgn}{=}&\!\!\!\frac{u_2^{-\theta}+\dots +u_n^{-\theta}-n+1}{u_1^{-\theta}+\dots +u_n^{-\theta}-n+1}<0,
		\end{eqnarray*}
		where the inequality follows from $\theta \in[-1,0)$. And for any $l\neq i$,
		\begin{eqnarray*}
			&&\!\!\!\frac{\partial \alpha_i(u_1,u_2,\dots,u_n)}{\partial u_l}\\
			&=&\!\!\!\theta u_i^{-2\theta-1}(u_1^{-\theta}+\dots +u_n^{-\theta}-n+1)^{-2}<0,~~\theta\in[-1,0),
		\end{eqnarray*}
		i.e., the function $\alpha_i(u_1,u_2,\dots,u_n)$ is decreasing in $\mathbf{u}\in [0,1]^n, ~i=1,2,\dots,n$.
		\item For the FGM family of Copulas,
		\begin{equation*}
			H(u_1,u_2,\dots, u_n)=\prod\limits_{k=1}^{n}u_k\left( 1+\theta\prod\limits_{k=1}^{n}(1-u_k)\right) ,~u_k\in[0,1].
		\end{equation*}
		Then
		\begin{eqnarray*}
			\alpha_i(u_1,u_2,\dots,u_n)&=&\frac{u_iD_{(i)}H(u_1,u_2,\dots,u_n)}{H(u_1,u_2,\dots,u_n)}\\
			&=&\frac{\left( \prod\limits_{k\neq i}^{n}u_k(1+\theta\prod\limits_{k=1}^{n}(1-u_k))-\theta\prod\limits_{k=1}^{n}u_k\prod\limits_{k\neq i}^{n}(1-u_k)\right) u_i }{\prod\limits_{k=1}^{n}u_k\left( 1+\theta\prod\limits_{k=1}^{n}(1-u_k)\right)}\\
			&=&1-\frac{\theta\prod\limits_{k\neq i}^{n}(1-u_k)u_i }{1+\theta\prod\limits_{k=1}^{n}(1-u_k)}.
		\end{eqnarray*}
		Taking the partial derivative of $\alpha_i(u_1,u_2,\dots,u_n)$ with respect to $u_i,$ we have
		\begin{eqnarray*}
			&&\!\!\!\frac{\partial \alpha_i(u_1,u_2,\dots,u_n)}{\partial u_i}\\
			&\overset{\rm sgn}{=}&\!\!\!-\theta\prod\limits_{k\neq i}^{n}(1-u_k)(1+\theta\prod\limits_{k=1}^{n}(1-u_k))-\theta u_i\prod\limits_{k\neq i}^{n}(1-u_k)\theta\prod\limits_{k\neq i}^{n}(1-u_k)\\
			&=&\!\!\!-\theta\prod\limits_{k\neq i}^{n}(1-u_k)\left[ 1+\theta\prod\limits_{k\neq i}^{n}(1-u_k)\right].
		\end{eqnarray*}
		Then $\alpha_i(u_1,u_2,\dots,u_n)$ decreases with respect to $u_i$ if and only if $\theta\in [0,1]$. And for any $l\neq i$,
		\begin{eqnarray*}
			&&\!\!\!\frac{\partial \alpha_i(u_1,u_2,\dots,u_n)}{\partial u_l}\\
			&\overset{\rm sgn}{=}&\!\!\!\theta u_i\prod\limits_{k\neq i,l}^{n}(1-u_k)\left[ 1+\theta\prod\limits_{k=1}^{n}(1-u_k) \right]-\theta^2 u_i\prod\limits_{k\neq i}^{n}(1-u_k)\prod\limits_{k\neq l}^{n}(1-u_k)\\
			&=&\theta u_i\prod\limits_{k\neq i,l}^{n}(1-u_k)\left[ 1+\theta\prod\limits_{k=1}^{n}(1-u_k)-(1-u_l)\prod\limits_{k\neq l}^{n}(1-u_k)\right]\\
			&=&\theta u_i\prod\limits_{k\neq i,l}^{n}(1-u_k)>0,~~ \theta\in[0,1].
		\end{eqnarray*}
		Therefore, the function $\alpha_i(u_1,u_2,\dots,u_n)$ is increasing in $\mathbf{u}\in[0,1]^n$ if and only if $\theta=0$.
	\end{enumerate}
\textbf{Proof of Proposition \ref{p2}.}
	Obviously, it is only necessary to prove the case where $n u^{-\theta}-n+1\geqslant0$, for the Clayton family of Copulas, by \cite{torrado2022optimal}, 
	\begin{eqnarray*}
		\eta(u)&=&\frac{n(1-u) u^{-\theta-1}\left(n u^{-\theta}-n+1\right)^{-\frac{1}{\theta}-1}}{1-\left(n u^{-\theta}-n+1\right)^{-\frac{1}{\theta}}}\\
		&=&\frac{n}{\left[ \left(n u^{-\theta}-n+1\right)^{\frac{1}{\theta}}-1\right]\left(n u^{-\theta}-n+1\right)u^{\theta+1}(1-u)^{-1} }.
	\end{eqnarray*}
	Differentiating the function defined $\eta(u)$ with respect to u, we get
	\begin{eqnarray*}
		\eta^{\prime}(u)\!\!\!\!&\overset{\rm sgn}{=}&\!\!\!\! n\theta\left( 1-u\right)\left( \left(n u^{-\theta}-n+1\right)^{\frac{1}{\theta}}-1\right)+n\left( 1-u\right) \left(n u^{-\theta}-n+1\right)^{\frac{1}{\theta}}-\left( \theta+1\right)u^\theta\\
		&&\!\!\!\!\times\left( 1-u\right)\left(n u^{-\theta}-n+1\right) \left( \left(n u^{-\theta}-n+1\right)^{\frac{1}{\theta}}-1\right)-u^{\theta+1}\left( 1-u\right)\left(n u^{-\theta}-n+1\right) \\
		&&\!\!\!\!\times\left( \left(n u^{-\theta}-n+1\right)^{\frac{1}{\theta}}-1\right)\\ 
		&{=}&\!\!\!\! n-n u\left(n u^{-\theta}-n+1\right)^{1 / \theta}+u^\theta(n-1)(1+\theta(1-u))
		\left(\left(n u^{-\theta}-n+1\right)^{1 / \theta}-1\right),
	\end{eqnarray*}
	where $\theta\in[-1,0)$. For convenience, we denote
	\begin{eqnarray*}
		\rho(u)&=&n-n u\left(n u^{-\theta}-n+1\right)^{1 / \theta}+u^\theta(n-1)(1+\theta(1-u))
		\left(\left(n u^{-\theta}-n+1\right)^{1 / \theta}-1\right).
	\end{eqnarray*}
	Note that $\rho(0)=n$ and $\rho(1)=0$. Furthermore,
	\begin{eqnarray*}
		\rho^{\prime}(u)&=&n(\theta-1)\left(n u^{-\theta}-n+1\right)^{1 / \theta}-n\theta-nu^\theta(n-1)\left[\left(n u^{-\theta}-n+1\right)^{1 / \theta}-1 \right]\\
		&=&n\theta\left[\left(n u^{-\theta}-n+1\right)^{1 / \theta}-1 \right]-n\left(n u^{-\theta}-n+1\right)^{1 / \theta}-nu^\theta(n-1)\\
		&&\times\left[\left(n u^{-\theta}-n+1\right)^{1 / \theta}-1 \right]\leqslant0, ~\text{for} ~\theta\in[-1,0).
	\end{eqnarray*}
	Then, the function $\rho(u)$ is monotone and $\rho(u)\geqslant0$. Hence, $\eta(u)$ is strictly increasing for $u\in(0,1)$ and $\theta\in[-1,0)$.
\appendixsection{Appendix B}
\label{appendix:B}
The proofs of the theorems in Section \ref{p} is provided here.\newline
\textbf{Proof of Theorem \ref{5}.} 
	For the inequality $C_s(K+1)-C_s(K)\geq0$,
	\begin{eqnarray}\label{(14)}
		&&\!\!\!\!\!\!\!\!\!\!\!\!\!\!\!\!\!\!\!\!\!\!\!\!\!\!\!\!\!\!\!\!\!\!\!\!\frac{c_f-(c_f-\sum\limits_{i=1}^{n}c_{p_i})\hat{C}(\overline{F}_1((K+1)\tau),\dots,\overline{F}_n((K+1)\tau))}{\int_{0}^{(K+1)\tau}\hat{C}(\overline{F}_1(t),\dots,\overline{F}_n(t))\mathrm{d}t}\nonumber\\
		&&\!\!\!\!\!\!\!\!\!\!\!\!\!\!\!\!\!\!\!\!\!\!\!\!\!\!\!\!\!\!\!\!\!\!\!\!\geq\frac{c_f-(c_f-\sum\limits_{i=1}^{n}c_{p_i})\hat{C}(\overline{F}_1(K\tau),\dots,\overline{F}_n(K\tau))}{\int_{0}^{K\tau}\hat{C}(\overline{F}_1(t),\dots,\overline{F}_n(t))\mathrm{d}t}\nonumber\\
		\text{if and only if}\nonumber\\
		&&\!\!\!\!\!\!\!\!\!\!\!\!\!\!\!\!\!\!\!\!\!\!\!\!\!\!\!\!\!\!\!\!\!\!\!\!N_s(K)\int_{0}^{K\tau}\hat{C}(\overline{F}_1(t),\dots,\overline{F}_n(t))\mathrm{d}t+\hat{C}(\overline{F}_1(K\tau),\dots,\overline{F}_n(K\tau))\geq\frac{c_f}{c_f-\sum\limits_{i=1}^{n}c_{p_i}},\nonumber\\
	\end{eqnarray}
	where 
	\begin{equation*}
		N_s(K)=\frac{\hat{C}(\overline{F}_1(K\tau),\dots,\overline{F}_n(K\tau))-\hat{C}(\overline{F}_1((K+1)\tau),\dots,\overline{F}_n((K+1)\tau))}{\int_{K\tau}^{(K+1)\tau}\hat{C}(\overline{F}_1(t),\dots,\overline{F}_n(t))\mathrm{d}t}.
	\end{equation*}
	It is known that $X_{1:n}$ is IFR by the proof of Theorem \ref{1}. Noting that
	\begin{equation*}
		h_{1:n}(K\tau)\leq N_s(K)\leq	h_{1:n}((K+1)\tau),
	\end{equation*}
	that is, $N_s(K)$ is strictly increasing with $K$.
	For convenience, let
	\begin{equation*}
		H_s(K)=N_s(K)\int_{0}^{K\tau}\hat{C}(\overline{F}_1(t),\dots,\overline{F}_n(t))\mathrm{d}t+\hat{C}(\overline{F}_1(K\tau),\dots,\overline{F}_n(K\tau)),
	\end{equation*}
	it should be show that $H_s(K)$ increases with $K$.
	For given $\tau>0$, 
	\begin{eqnarray*}
		&&H_s(K+1 )-H_s(K)\\
		&=&N_s(K+1)\int_{0}^{(K+1)\tau}\hat{C}(\overline{F}_1(t),\dots,\overline{F}_n(t))\mathrm{d}t+\hat{C}(\overline{F}_1((K+1)\tau),\dots,\overline{F}_n((K+1)\tau))\nonumber\\
		&&-N_s(K)\int_{0}^{K\tau}\hat{C}(\overline{F}_1(t),\dots,\overline{F}_n(t))\mathrm{d}t-\hat{C}(\overline{F}_1(K\tau),\dots,\overline{F}_n(K\tau))\\
		&=&\left[ N_s(K+1)-N_s(K)\right]\int_{0}^{(K+1)\tau}\hat{C}(\overline{F}_1(t),\dots,\overline{F}_n(t))\mathrm{d}t >0.
	\end{eqnarray*}
	Additionly, 
	\begin{eqnarray*}
		H_s(\infty)=\lim\limits_{K\to \infty}H_s(K)=h_s(\infty)\mu_s.
	\end{eqnarray*}
	Thus, there exists a finite and unique minimum $K^*~(1\leq K^*<\infty)$ which satisfies (\ref{(14)}).
\textbf{Proof of Theorem \ref{6}.} 
	For the inequality $C_{s1}(K+1)-C_{s1}(K)\geq0$,
	\begin{eqnarray}\label{(16)}
		&&\!\!\!\!\!\!\!\!\!\!\!\!\!\!\!\!\!\!\!\!\!\!\!\!\!\!\!\!\!\!\!\!\!\!\!\!\!\!\!\!\!\!\!\!\!\!\!\!\!\!\frac{c_f\!\!-\!\!(c_f\!-\!\sum\limits_{i=1}^{n}c_{p_i})\hat{C}(\overline{F}_1((K\!+\!1)\tau),\!\dots,\!\overline{F}_n((K\!+\!1)\tau))\!+\!	c_{d1}  \int_{0}^{(K+1)\tau}( 1\!-\!\hat{C}(\overline{F}_1(t),\!\dots,\!\overline{F}_n(t) \mathrm{d}t))\!+\!c_{d2}\mu_s}{\int_{0}^{(K+1)\tau}\hat{C}(\overline{F}_1(t),\!\dots,\!\overline{F}_n(t))\mathrm{d}t}\nonumber\\
		&&\!\!\!\!\!\!\!\!\!\!\!\!\!\!\!\!\!\!\!\!\!\!\!\!\!\!\!\!\!\!\!\!\!\!\!\!\!\!\!\!\!\!\!\!\!\!\!\!\!\!\geq\frac{c_f\!-\!(c_f\!-\!\sum\limits_{i=1}^{n}c_{p_i})\hat{C}(\overline{F}_1(K\tau),\!\dots,\!\overline{F}_n(K\tau))\!+\!c_{d1}  \int_{0}^{K\tau}( 1\!-\!\hat{C}(\overline{F}_1(t),\!\dots,\!\overline{F}_n(t))) \mathrm{d}t\!+\!c_{d2}\mu_s}{\int_{0}^{K\tau}\hat{C}(\overline{F}_1(t),\!\dots,\!\overline{F}_n(t))\mathrm{d}t}\nonumber\\
		\text{if and only if}\nonumber\\
		&&\!\!\!\!\!\!\!\!\!\!\!\!\!\!\!\!\!\!\!\!\!\!\!\!\!\!\!\!\!\!\!\!\!\!\!\!\!\!\!(c_f-\sum\limits_{i=1}^{n}c_{p_i})H_s(K)+c_{d1} J_s(K) \geq c_f+c_{d2}\mu_s,
	\end{eqnarray}
	where 
	\begin{eqnarray*}
		&&J_s(K)=M_s(K)\int_{0}^{K\tau}\hat{C}(\overline{F}_1(t),\dots,\overline{F}_n(t)) \mathrm{d}t-\int_{0}^{K\tau}( 1-\hat{C}(\overline{F}_1(t),\!\dots,\!\overline{F}_n(t))) \mathrm{d}t,\\
		&&M_s(K)=\frac{\int_{K\tau}^{(K+1)\tau}( 1-\hat{C}(\overline{F}_1(t),\dots,\overline{F}_n(t))) \mathrm{d}t}{\int_{K\tau}^{(K+1)\tau}\hat{C}(\overline{F}_1(t),\dots,\overline{F}_n(t)) \mathrm{d}t}.
	\end{eqnarray*}
	Noting that 
	\begin{equation*}
		\frac{1-\hat{C}(\overline{F}_1(K\tau),\dots,\overline{F}_n(K\tau))}{\hat{C}(\overline{F}_1(K\tau),\dots,\overline{F}_n(K\tau))}<M_s(K)<\frac{1-\hat{C}(\overline{F}_1((K+1)\tau),\!\dots,\!\overline{F}_n((K+1)\tau))}{\hat{C}(\overline{F}_1((K+1)\tau),\!\dots,\!\overline{F}_n((K+1)\tau))},
	\end{equation*}
	that is, $M_s(K)$ is strictly increasing with $K$. For given $\tau>0$,
	\begin{eqnarray*}
		&&\!\!\!J_s(K+1)-J_s(K)\\
		&=&\!\!\!M_s(K+1)\int_{0}^{(K+1)\tau}\hat{C}(\overline{F}_1(t),\dots,\overline{F}_n(t)) \mathrm{d}t-M_s(K)\int_{0}^{K\tau}\hat{C}(\overline{F}_1(t),\dots,\overline{F}_n(t)) \mathrm{d}t\\
		&&-\int_{K\tau}^{(K+1)\tau}( 1-\hat{C}(\overline{F}_1(t),\dots,\overline{F}_n(t))) \mathrm{d}t\\
		&>&M_s(K)\int_{K\tau}^{(K+1)\tau}\hat{C}(\overline{F}_1(t),\dots,\overline{F}_n(t)) \mathrm{d}t-\int_{K\tau}^{(K+1)\tau}( 1-\hat{C}(\overline{F}_1(t),\dots,\overline{F}_n(t))) \mathrm{d}t=0.
	\end{eqnarray*}
	Hence, the left-hand side of (\ref{(16)}) increases strictly with $K$ to $\infty$. 	Thus, there exists a finite and unique minimum $K^*~(1\leq K^*<\infty)$ which satisfies (\ref{(16)}).\newline
\textbf{Proof of Theorem \ref{7}.} 
	For the inequality $C_p(K+1)-C_p(K)\geq0$,
	\begin{eqnarray}\label{(18)}
		&&\!\!\!\!\!\!\!\!\!\!\!\!\!\!\!\!\!\!\!\!\!\!\!\!\!\!\!\!\!\!\!\!\!\!\!\!\frac{\sum\limits_{i=1}^{n}c_{p_i}+(c_f-\sum\limits_{i=1}^{n}c_{p_i}){C}({F}_1((K+1)\tau),\dots,{F}_n((K+1)\tau))}{\int_{0}^{(K+1)\tau}(1-{C}({F}_1(t),\dots,{F}_n(t)))\mathrm{d}t}\nonumber\\
		&&\!\!\!\!\!\!\!\!\!\!\!\!\!\!\!\!\!\!\!\!\!\!\!\!\!\!\!\!\!\!\!\!\!\!\!\!\geq\frac{\sum\limits_{i=1}^{n}c_{p_i}+(c_f-\sum\limits_{i=1}^{n}c_{p_i}){C}({F}_1(K\tau),\dots,{F}_n(K\tau))}{\int_{0}^{K\tau}(1-{C}({F}_1(t),\dots,{F}_n(t)))\mathrm{d}t}\nonumber\\
		\text{if and only if}\nonumber\\
		&&\!\!\!\!\!\!\!\!\!\!\!\!\!\!\!\!\!\!\!\!\!\!\!\!\!\!\!\!\!\!\!\!\!\!\!\!H_p(K)\geq\frac{\sum\limits_{i=1}^{n}c_{p_i}}{c_f-\sum\limits_{i=1}^{n}c_{p_i}},
	\end{eqnarray}
	where 
	\begin{eqnarray*}
		\!\!\!\!&&\!\!\!\!H_p(K)=N_p(K)\int_{0}^{K\tau}(1-{C}({F}_1(t),\dots,{F}_n(t)))\mathrm{d}t-{C}({F}_1(K\tau),\dots,{F}_n(K\tau)),\\
		\!\!\!\!&&\!\!\!\!N_p(K)=\frac{{C}({F}_1((K+1)\tau),\dots,{F}_n((K+1)\tau))-{C}({F}_1(K\tau),\dots,{F}_n(K\tau))}{\int_{K\tau}^{(K+1)\tau}(1-{C}({F}_1(t),\dots,{F}_n(t)))\mathrm{d}t}.
	\end{eqnarray*}
	It is known that $X_{n:n}$ is IFR by the proof of Theorem \ref{3}. Noting that
	\begin{equation*}
		h_{n:n}(K\tau)\leq N_p(K)\leq	h_{n:n}((K+1)\tau),
	\end{equation*}
	that is, $N_p(K)$ is strictly increasing with $K$. For given $\tau>0$,
	\begin{eqnarray*}
		&&H_p(K+1 )-H_p(K)\\
		&=&N_p(K)\int_{0}^{(K+1)\tau}(1-{C}({F}_1(t),\dots,{F}_n(t)))\mathrm{d}t-{C}({F}_1((K+1)\tau),\dots,{F}_n((K+1)\tau))\nonumber\\
		&&-N_p(K)\int_{0}^{K\tau}(1-{C}({F}_1(t),\dots,{F}_n(t)))\mathrm{d}t+{C}({F}_1(K\tau),\dots,{F}_n(K\tau))\nonumber\\
		&=&\left[ N_p(K+1)-N_p(K)\right]\int_{0}^{(K+1)\tau}(1-{C}({F}_1(t),\dots,{F}_n(t)))\mathrm{d}t \nonumber\\
		&>&0.
	\end{eqnarray*}
	Additionly, 
	\begin{eqnarray*}
		H_p(\infty)=\lim\limits_{K\to \infty}H_p(K)=h_p(\infty)\mu_p-1.
	\end{eqnarray*}
	Thus, there exists a finite and unique minimum $K^*~(1\leq K^*<\infty)$ which satisfies (\ref{(18)}) if $h_{n:n}(\infty)\mu_p>c_f/(c_f-\sum_{i=1}^{n}c_{p_i})$.\newline
\textbf{Proof of Theorem \ref{8}.} 
	For the inequality $C_{p1}(K+1)-C_{p1}(K)\geq0$,
	\begin{eqnarray}\label{(20)}
		&&\!\!\!\!\!\!\!\!\!\!\!\!\!\!\!\!\!\!\!\!\!\!\!\!\!\!\!\!\!\!\!\!\!\!\!\!\!\!\!\!\!\!\!\!\!\!\!\!\!\!\!\frac{\sum\limits_{i=1}^{n}c_{p_i}\!\!+\!\!(c_f\!-\!\sum\limits_{i=1}^{n}\!c_{p_i}){C}({F}_1((K\!+\!1)\tau),\!\dots\!,\!{F}_n((K\!+\!1)\tau))\!+\!	c_{d1}  \int_{0}^{(K+1)\tau}\!{C}({F}_1(t),\!\dots\!,\!{F}_n(t)) \mathrm{d}t\!+\!c_{d2}\mu_p}{\int_{0}^{(K+1)\tau}(1-{C}({F}_1(t),\!\dots\!,{F}_n(t)))\mathrm{d}t}\nonumber\\
		&&\!\!\!\!\!\!\!\!\!\!\!\!\!\!\!\!\!\!\!\!\!\!\!\!\!\!\!\!\!\!\!\!\!\!\!\!\!\!\!\!\!\!\!\!\!\!\!\!\!\!\!\geq\frac{\sum\limits_{i=1}^{n}c_{p_i}\!+\!(c_f\! -\!\sum\limits_{i=1}^{n}c_{p_i}){C}({F}_1(K\tau),\dots,{F}_n(K\tau))\!+\!	c_{d1}  \int_{0}^{K\tau}{C}({F}_1(t),\dots,{F}_n(t)) \mathrm{d}t\!+\!c_{d2}\mu_p}{\int_{0}^{K\tau}(1-{C}({F}_1(t),\dots,{F}_n(t)))\mathrm{d}t}\nonumber\\
		\text{is equivalent to}\nonumber\\
		&&\!\!\!\!\!\!\!\!\!\!\!\!\!\!\!\!\!\!\!\!\!\!\!\!\!\!\!\!\!\!\!\!\!\!\!\!\!\!\!\!\!\!\!\!\!\!\!\!\!\!\!(c_f-\sum\limits_{i=1}^{n}c_{p_i})H_p(K)+c_{d1} J_p(K) \geq \sum\limits_{i=1}^{n}c_{p_i}+c_{d2}\mu_p,
	\end{eqnarray}
	where 
	\begin{eqnarray*}
		&&J_p(K)=M_p(K)\int_{0}^{K\tau}(1-{C}({F}_1(t),\dots,{F}_n(t))) \mathrm{d}t-\int_{0}^{K\tau}{C}({F}_1(t),\dots,{F}_n(t)) \mathrm{d}t,\\
		&&M_p(K)=\frac{\int_{K\tau}^{(K+1)\tau}{C}({F}_1(t),\dots,{F}_n(t)) \mathrm{d}t}{\int_{K\tau}^{(K+1)\tau}(1-{C}({F}_1(t),\dots,{F}_n(t))) \mathrm{d}t}.
	\end{eqnarray*}
	Noting that 
	\begin{equation*}
		\frac{{C}({F}_1(K\tau),\dots,{F}_n(K\tau))}{1-{C}({F}_1(K\tau),\dots,{F}_n(K\tau))}<M_p(K)<\frac{{C}({F}_1((K+1)\tau),\dots,{F}_n((K+1)\tau))}{1-{C}({F}_1((K+1)\tau),\dots,{F}_n((K+1)\tau))},
	\end{equation*}
	that is, $M_p(K)$ is strictly increasing with $K$. For given $\tau>0$,
	\begin{eqnarray*}
		&&\!\!\!J_p(K+1)-J_p(K)\\
		&=&\!\!\!M_p(K+1)\!\int_{0}^{(K+1)\tau}\!(1\!-\!{C}({F}_1(t),\dots,{F}_n(t))) \mathrm{d}t\!-\!M_p(K)\!\int_{0}^{K\tau}\!(1\!-\!{C}({F}_1(t),\dots,{F}_n(t))) \mathrm{d}t\\
		&&\!\!\!-\int_{K\tau}^{(K+1)\tau}{C}({F}_1(t),\dots,{F}_n(t)) \mathrm{d}t\\
		&>&\!\!\!M_p(K)\int_{K\tau}^{(K+1)\tau}(1-{C}({F}_1(t),\dots,{F}_n(t))) \mathrm{d}t-\int_{K\tau}^{(K+1)\tau}{C}({F}_1(t),\dots,{F}_n(t)) \mathrm{d}t=0.
	\end{eqnarray*}
	Hence, the left-hand side of (\ref{(20)}) increases strictly with $K$ to $\infty$. 	Thus, there exists a finite and unique minimum $K^*~(1\leq K^*<\infty)$ which satisfies (\ref{(20)}).
\bibliographystyle{mystyle}
\bibliography{ref-quan}

\begin{thebibliography}{39}
\expandafter\ifx\csname natexlab\endcsname\relax\def\natexlab#1{#1}\fi
\providecommand{\url}[1]{\texttt{#1}}
\providecommand{\href}[2]{#2}
\providecommand{\path}[1]{#1}
\providecommand{\DOIprefix}{doi:}
\providecommand{\ArXivprefix}{arXiv:}
\providecommand{\URLprefix}{URL: }
\providecommand{\Pubmedprefix}{pmid:}
\providecommand{\doi}[1]{\href{http://dx.doi.org/#1}{\path{#1}}}
\providecommand{\Pubmed}[1]{\href{pmid:#1}{\path{#1}}}
\providecommand{\bibinfo}[2]{#2}
\ifx\xfnm\relax \def\xfnm[#1]{\unskip,\space#1}\fi
\bibitem[{Akhtar et~al.(2021)Akhtar, Kirmani \& Jameel}]{akhtar2021reliability}
\bibinfo{author}{Akhtar, I.}, \bibinfo{author}{Kirmani, S.}, \&
  \bibinfo{author}{Jameel, M.} (\bibinfo{year}{2021}).
\newblock \bibinfo{title}{Reliability assessment of power system considering
  the impact of renewable energy sources integration into grid with advanced
  intelligent strategies}.
\newblock {\it \bibinfo{journal}{IEEE Access}\/},  {\it \bibinfo{volume}{9}\/},
  \bibinfo{pages}{32485--32497}.
\bibitem[{Barlow \& Proschan(1996)}]{barlow1996mathematical}
\bibinfo{author}{Barlow, R.~E.}, \& \bibinfo{author}{Proschan, F.}
  (\bibinfo{year}{1996}).
\newblock {\it \bibinfo{title}{Mathematical theory of reliability}\/}.
\newblock \bibinfo{publisher}{SIAM}.
\bibitem[{Barroso \& Clidaras(2022)}]{barroso2022datacenter}
\bibinfo{author}{Barroso, L.~A.}, \& \bibinfo{author}{Clidaras, J.}
  (\bibinfo{year}{2022}).
\newblock {\it \bibinfo{title}{The datacenter as a computer: An introduction to
  the design of warehouse-scale machines}\/}.
\newblock \bibinfo{publisher}{Springer Nature}.
\bibitem[{Belzunce et~al.(2001)Belzunce, Franco, Ruiz \&
  Ruiz}]{belzunce2001partial}
\bibinfo{author}{Belzunce, F.}, \bibinfo{author}{Franco, M.},
  \bibinfo{author}{Ruiz, J.-M.}, \& \bibinfo{author}{Ruiz, M.~C.}
  (\bibinfo{year}{2001}).
\newblock \bibinfo{title}{On partial orderings between coherent systems with
  different structures}.
\newblock {\it \bibinfo{journal}{Probability in the Engineering and
  Informational Sciences}\/},  {\it \bibinfo{volume}{15}\/},
  \bibinfo{pages}{273--293}.
\bibitem[{Berg(1976)}]{berg1976proof}
\bibinfo{author}{Berg, M.} (\bibinfo{year}{1976}).
\newblock \bibinfo{title}{A proof of optimality for age replacement policies}.
\newblock {\it \bibinfo{journal}{Journal of Applied Probability}\/},  {\it
  \bibinfo{volume}{13}\/}, \bibinfo{pages}{751--759}.
\bibitem[{Eryilmaz(2023)}]{eryilmaz2023age}
\bibinfo{author}{Eryilmaz, S.} (\bibinfo{year}{2023}).
\newblock \bibinfo{title}{Age based preventive replacement policy for discrete
  time coherent systems with independent and identical components}.
\newblock {\it \bibinfo{journal}{Reliability Engineering \& System Safety}\/},
  {\it \bibinfo{volume}{240}\/}, \bibinfo{pages}{109544}.
\bibitem[{Eryilmaz \& Ozkut(2020)}]{eryilmaz2020optimization}
\bibinfo{author}{Eryilmaz, S.}, \& \bibinfo{author}{Ozkut, M.}
  (\bibinfo{year}{2020}).
\newblock \bibinfo{title}{Optimization problems for a parallel system with
  multiple types of dependent components}.
\newblock {\it \bibinfo{journal}{Reliability Engineering \& System Safety}\/},
  {\it \bibinfo{volume}{199}\/}, \bibinfo{pages}{106911}.
\bibitem[{Eryilmaz \& Tank(2023)}]{eryilmaz2023optimal}
\bibinfo{author}{Eryilmaz, S.}, \& \bibinfo{author}{Tank, F.}
  (\bibinfo{year}{2023}).
\newblock \bibinfo{title}{Optimal age replacement policy for discrete time
  parallel systems}.
\newblock {\it \bibinfo{journal}{Top}\/},  {\it \bibinfo{volume}{31}\/},
  \bibinfo{pages}{475--490}.
\bibitem[{Kundur(2007)}]{kundur2007power}
\bibinfo{author}{Kundur, P.} (\bibinfo{year}{2007}).
\newblock \bibinfo{title}{Power system stability}.
\newblock {\it \bibinfo{journal}{Power system stability and control}\/},  {\it
  \bibinfo{volume}{10}\/}, \bibinfo{pages}{7--1}.
\bibitem[{Levitin et~al.(2021)Levitin, Xing \& Dai}]{levitin2021optimal}
\bibinfo{author}{Levitin, G.}, \bibinfo{author}{Xing, L.}, \&
  \bibinfo{author}{Dai, Y.} (\bibinfo{year}{2021}).
\newblock \bibinfo{title}{Optimal operation and maintenance scheduling in
  m-out-of-n standby systems with reusable elements}.
\newblock {\it \bibinfo{journal}{Reliability Engineering \& System Safety}\/},
  {\it \bibinfo{volume}{211}\/}, \bibinfo{pages}{107582}.
\bibitem[{Levitin et~al.(2023)Levitin, Xing \& Dai}]{levitin2023standby}
\bibinfo{author}{Levitin, G.}, \bibinfo{author}{Xing, L.}, \&
  \bibinfo{author}{Dai, Y.} (\bibinfo{year}{2023}).
\newblock \bibinfo{title}{Standby mode transfer schedule minimizing downtime of
  1-out-of-n system with storage}.
\newblock {\it \bibinfo{journal}{Reliability Engineering \& System Safety}\/},
  {\it \bibinfo{volume}{237}\/}, \bibinfo{pages}{109322}.
\bibitem[{Levitin et~al.(2024)Levitin, Xing \& Dai}]{levitin2024optimizing}
\bibinfo{author}{Levitin, G.}, \bibinfo{author}{Xing, L.}, \&
  \bibinfo{author}{Dai, Y.} (\bibinfo{year}{2024}).
\newblock \bibinfo{title}{Optimizing corrective maintenance for multistate
  systems with storage}.
\newblock {\it \bibinfo{journal}{Reliability Engineering \& System Safety}\/},
  {\it \bibinfo{volume}{244}\/}, \bibinfo{pages}{109951}.
\bibitem[{Li \& Wu(2024)}]{li2024optimal}
\bibinfo{author}{Li, M.}, \& \bibinfo{author}{Wu, B.} (\bibinfo{year}{2024}).
\newblock \bibinfo{title}{Optimal condition-based opportunistic maintenance
  policy for two-component systems considering common cause failure}.
\newblock {\it \bibinfo{journal}{Reliability Engineering \& System Safety}\/},
  (p. \bibinfo{pages}{110269}).
\bibitem[{Liu \& Wang(2021)}]{liu2021optimal}
\bibinfo{author}{Liu, P.}, \& \bibinfo{author}{Wang, G.}
  (\bibinfo{year}{2021}).
\newblock \bibinfo{title}{Optimal periodic preventive maintenance policies for
  systems subject to shocks}.
\newblock {\it \bibinfo{journal}{Applied Mathematical Modelling}\/},  {\it
  \bibinfo{volume}{93}\/}, \bibinfo{pages}{101--114}.
\bibitem[{Medara \& Singh(2021)}]{medara2021energy}
\bibinfo{author}{Medara, R.}, \& \bibinfo{author}{Singh, R.~S.}
  (\bibinfo{year}{2021}).
\newblock \bibinfo{title}{Energy efficient and reliability aware workflow task
  scheduling in cloud environment}.
\newblock {\it \bibinfo{journal}{Wireless Personal Communications}\/},  {\it
  \bibinfo{volume}{119}\/}, \bibinfo{pages}{1301--1320}.
\bibitem[{Nakagawa(2006)}]{nakagawa2006maintenance}
\bibinfo{author}{Nakagawa, T.} (\bibinfo{year}{2006}).
\newblock {\it \bibinfo{title}{Maintenance theory of reliability}\/}.
\newblock \bibinfo{publisher}{Springer Science \& Business Media}.
\bibitem[{Nakagawa(2008)}]{nakagawa2008advanced}
\bibinfo{author}{Nakagawa, T.} (\bibinfo{year}{2008}).
\newblock {\it \bibinfo{title}{Advanced reliability models and maintenance
  policies}\/}.
\newblock \bibinfo{publisher}{Springer Science \& Business Media}.
\bibitem[{Nakagawa \& Zhao(2012)}]{nakagawa2012optimization}
\bibinfo{author}{Nakagawa, T.}, \& \bibinfo{author}{Zhao, X.}
  (\bibinfo{year}{2012}).
\newblock \bibinfo{title}{Optimization problems of a parallel system with a
  random number of units}.
\newblock {\it \bibinfo{journal}{IEEE Transactions on Reliability}\/},  {\it
  \bibinfo{volume}{61}\/}, \bibinfo{pages}{543--548}.
\bibitem[{Navarro et~al.(2014)Navarro, del {\'A}guila, Sordo \&
  Su{\'a}rez-Llorens}]{navarro2014preservation}
\bibinfo{author}{Navarro, J.}, \bibinfo{author}{del {\'A}guila, Y.},
  \bibinfo{author}{Sordo, M.~A.}, \& \bibinfo{author}{Su{\'a}rez-Llorens, A.}
  (\bibinfo{year}{2014}).
\newblock \bibinfo{title}{Preservation of reliability classes under the
  formation of coherent systems}.
\newblock {\it \bibinfo{journal}{Applied Stochastic Models in Business and
  Industry}\/},  {\it \bibinfo{volume}{30}\/}, \bibinfo{pages}{444--454}.
\bibitem[{Nelsen(2006)}]{Nelsen2006}
\bibinfo{author}{Nelsen, R.~B.} (\bibinfo{year}{2006}).
\newblock {\it \bibinfo{title}{An Introduction to Copulas (Springer Series in
  Statistics)}\/} volume~\bibinfo{volume}{47}.
\newblock \bibinfo{publisher}{Springer-Verlag Berlin, Heidelberg}.
\bibitem[{Ota \& Kimura(2017)}]{ota2017statistical}
\bibinfo{author}{Ota, S.}, \& \bibinfo{author}{Kimura, M.}
  (\bibinfo{year}{2017}).
\newblock \bibinfo{title}{A statistical dependent failure detection method for
  n-component parallel systems}.
\newblock {\it \bibinfo{journal}{Reliability Engineering \& System Safety}\/},
  {\it \bibinfo{volume}{167}\/}, \bibinfo{pages}{376--382}.
\bibitem[{Safaei et~al.(2020)Safaei, Ch{\^a}telet \&
  Ahmadi}]{safaei2020optimal}
\bibinfo{author}{Safaei, F.}, \bibinfo{author}{Ch{\^a}telet, E.}, \&
  \bibinfo{author}{Ahmadi, J.} (\bibinfo{year}{2020}).
\newblock \bibinfo{title}{Optimal age replacement policy for parallel and
  series systems with dependent components}.
\newblock {\it \bibinfo{journal}{Reliability Engineering \& System Safety}\/},
  {\it \bibinfo{volume}{197}\/}, \bibinfo{pages}{106798}.
\bibitem[{Sheu et~al.(2018)Sheu, Liu, Zhang \& Tsai}]{sheu2018generalized}
\bibinfo{author}{Sheu, S.-H.}, \bibinfo{author}{Liu, T.-H.},
  \bibinfo{author}{Zhang, Z.-G.}, \& \bibinfo{author}{Tsai, H.-N.}
  (\bibinfo{year}{2018}).
\newblock \bibinfo{title}{The generalized age maintenance policies with random
  working times}.
\newblock {\it \bibinfo{journal}{Reliability Engineering \& System Safety}\/},
  {\it \bibinfo{volume}{169}\/}, \bibinfo{pages}{503--514}.
\bibitem[{Torrado(2022)}]{torrado2022optimal}
\bibinfo{author}{Torrado, N.} (\bibinfo{year}{2022}).
\newblock \bibinfo{title}{Optimal component-type allocation and replacement
  time policies for parallel systems having multi-types dependent components}.
\newblock {\it \bibinfo{journal}{Reliability Engineering \& System Safety}\/},
  {\it \bibinfo{volume}{224}\/}, \bibinfo{pages}{108502}.
\bibitem[{Wang et~al.(2024)Wang, Wang, Zhao \& Miao}]{wang2024optimization}
\bibinfo{author}{Wang, J.}, \bibinfo{author}{Wang, L.}, \bibinfo{author}{Zhao,
  X.}, \& \bibinfo{author}{Miao, Z.} (\bibinfo{year}{2024}).
\newblock \bibinfo{title}{Optimization problems and maintenance policy for a
  parallel computing system with dependent components}.
\newblock {\it \bibinfo{journal}{Annals of Operations Research}\/},  (pp.
  \bibinfo{pages}{1--26}).
\bibitem[{Wang et~al.(2022)Wang, Ye \& Wang}]{wang2022extended}
\bibinfo{author}{Wang, J.}, \bibinfo{author}{Ye, J.}, \& \bibinfo{author}{Wang,
  L.} (\bibinfo{year}{2022}).
\newblock \bibinfo{title}{Extended age maintenance models and its optimization
  for series and parallel systems}.
\newblock {\it \bibinfo{journal}{Annals of Operations Research}\/},  (pp.
  \bibinfo{pages}{1--23}).
\bibitem[{Wireman(2004)}]{wireman2004total}
\bibinfo{author}{Wireman, T.} (\bibinfo{year}{2004}).
\newblock {\it \bibinfo{title}{Total productive maintenance}\/}.
\newblock \bibinfo{publisher}{Industrial Press Inc.}
\bibitem[{Wu \& Scarf(2017)}]{wu2017two}
\bibinfo{author}{Wu, S.}, \& \bibinfo{author}{Scarf, P.}
  (\bibinfo{year}{2017}).
\newblock \bibinfo{title}{Two new stochastic models of the failure process of a
  series system}.
\newblock {\it \bibinfo{journal}{European Journal of Operational Research}\/},
  {\it \bibinfo{volume}{257}\/}, \bibinfo{pages}{763--772}.
\bibitem[{Xing(2024)}]{xing2024decision}
\bibinfo{author}{Xing, L.} (\bibinfo{year}{2024}).
\newblock \bibinfo{title}{Decision diagrams for complex system reliability
  analysis}.
\newblock In {\it \bibinfo{booktitle}{Frontiers of Performability Engineering:
  In Honor of Prof. KB Misra}\/} (pp. \bibinfo{pages}{51--67}).
\newblock \bibinfo{publisher}{Springer}.
\bibitem[{Xing et~al.(2020)Xing, Zhao, Wang \& Xiang}]{xing2020reliability}
\bibinfo{author}{Xing, L.}, \bibinfo{author}{Zhao, G.}, \bibinfo{author}{Wang,
  Y.}, \& \bibinfo{author}{Xiang, Y.} (\bibinfo{year}{2020}).
\newblock \bibinfo{title}{Reliability modeling of correlated competitions and
  dependent components with random failure propagation time}.
\newblock {\it \bibinfo{journal}{Quality and Reliability Engineering
  International}\/},  {\it \bibinfo{volume}{36}\/}, \bibinfo{pages}{947--964}.
\bibitem[{Xing et~al.(2021)Xing, Zhao, Xiang \& Liu}]{xing2021behavior}
\bibinfo{author}{Xing, L.}, \bibinfo{author}{Zhao, G.}, \bibinfo{author}{Xiang,
  Y.}, \& \bibinfo{author}{Liu, Q.} (\bibinfo{year}{2021}).
\newblock \bibinfo{title}{A behavior-driven reliability modeling method for
  complex smart systems}.
\newblock {\it \bibinfo{journal}{Quality and Reliability Engineering
  International}\/},  {\it \bibinfo{volume}{37}\/},
  \bibinfo{pages}{2065--2084}.
\bibitem[{Yan \& Wang(2022)}]{yan2022component}
\bibinfo{author}{Yan, R.}, \& \bibinfo{author}{Wang, J.}
  (\bibinfo{year}{2022}).
\newblock \bibinfo{title}{Component level versus system level at active
  redundancies for coherent systems with dependent heterogeneous components}.
\newblock {\it \bibinfo{journal}{Communications in Statistics-Theory and
  Methods}\/},  {\it \bibinfo{volume}{51}\/}, \bibinfo{pages}{1724--1744}.
\bibitem[{Zeng et~al.(2023)Zeng, Barros \& Coit}]{zeng2023dependent}
\bibinfo{author}{Zeng, Z.}, \bibinfo{author}{Barros, A.}, \&
  \bibinfo{author}{Coit, D.} (\bibinfo{year}{2023}).
\newblock \bibinfo{title}{Dependent failure behavior modeling for risk and
  reliability: A systematic and critical literature review}.
\newblock {\it \bibinfo{journal}{Reliability Engineering \& System Safety}\/},
  (p. \bibinfo{pages}{109515}).
\bibitem[{Zhang et~al.(2023)Zhang, Feng \& Chen}]{zhang2023preventive}
\bibinfo{author}{Zhang, J.}, \bibinfo{author}{Feng, H.}, \&
  \bibinfo{author}{Chen, X.} (\bibinfo{year}{2023}).
\newblock \bibinfo{title}{Preventive maintenance policies for a big data system
  with throughput rate}.
\newblock {\it \bibinfo{journal}{Annals of Operations Research}\/},  (pp.
  \bibinfo{pages}{1--24}).
\bibitem[{Zhao et~al.(2024)Zhao, Bu, Pang \& Cai}]{zhao2024periodic}
\bibinfo{author}{Zhao, X.}, \bibinfo{author}{Bu, Y.}, \bibinfo{author}{Pang,
  W.}, \& \bibinfo{author}{Cai, J.} (\bibinfo{year}{2024}).
\newblock \bibinfo{title}{Periodic and random incremental backup policies in
  reliability theory}.
\newblock {\it \bibinfo{journal}{Software Quality Journal}\/},  (pp.
  \bibinfo{pages}{1--16}).
\bibitem[{Zhao et~al.(2022)Zhao, Mizutani, Chen \&
  Nakagawa}]{zhao2022preventive}
\bibinfo{author}{Zhao, X.}, \bibinfo{author}{Mizutani, S.},
  \bibinfo{author}{Chen, M.}, \& \bibinfo{author}{Nakagawa, T.}
  (\bibinfo{year}{2022}).
\newblock \bibinfo{title}{Preventive replacement policies for parallel systems
  with deviation costs between replacement and failure}.
\newblock {\it \bibinfo{journal}{Annals of Operations Research}\/},  (pp.
  \bibinfo{pages}{1--19}).
\bibitem[{Zhao et~al.(2015)Zhao, Mizutani \& Nakagawa}]{zhao2015better}
\bibinfo{author}{Zhao, X.}, \bibinfo{author}{Mizutani, S.}, \&
  \bibinfo{author}{Nakagawa, T.} (\bibinfo{year}{2015}).
\newblock \bibinfo{title}{Which is better for replacement policies with
  continuous or discrete scheduled times?}
\newblock {\it \bibinfo{journal}{European Journal of Operational Research}\/},
  {\it \bibinfo{volume}{242}\/}, \bibinfo{pages}{477--486}.
\bibitem[{Zhao \& Nakagawa(2012)}]{zhao2012optimization}
\bibinfo{author}{Zhao, X.}, \& \bibinfo{author}{Nakagawa, T.}
  (\bibinfo{year}{2012}).
\newblock \bibinfo{title}{Optimization problems of replacement first or last in
  reliability theory}.
\newblock {\it \bibinfo{journal}{European journal of operational research}\/},
  {\it \bibinfo{volume}{223}\/}, \bibinfo{pages}{141--149}.
\bibitem[{Zhao et~al.(2014)Zhao, Nakagawa \& Zuo}]{zhao2014optimal}
\bibinfo{author}{Zhao, X.}, \bibinfo{author}{Nakagawa, T.}, \&
  \bibinfo{author}{Zuo, M.~J.} (\bibinfo{year}{2014}).
\newblock \bibinfo{title}{Optimal replacement last with continuous and discrete
  policies}.
\newblock {\it \bibinfo{journal}{IEEE Transactions on Reliability}\/},  {\it
  \bibinfo{volume}{63}\/}, \bibinfo{pages}{868--880}.

\end{thebibliography}

\end{document}